\documentclass[a4paper,twoside]{article}
\usepackage[latin1]{inputenc} 
\usepackage[T1]{fontenc} 
\usepackage{RR}

\usepackage{amsmath}
\usepackage{amssymb}
\usepackage{graphicx}

\usepackage{moreverb}

\usepackage{subfig} 
\usepackage{multirow} 

\newcommand\BibTeX{{\rmfamily B\kern-.05em \textsc{i\kern-.025em b}\kern-.08em
T\kern-.1667em\lower.7ex\hbox{E}\kern-.125emX}}

\usepackage[normalem]{ulem}

\RRdate{November 2011}
\RRauthor{
Muriel Boulakia
\and
Elisa Schenone
\and
Jean-Fr\'ederic Gerbeau}
\authorhead{M. Boulakia \& E. Schenone \& J-F. Gerbeau}
\RRtitle{Mod\`eles r\'eduits pour l'\'electrophysiologie cardiaque. Application \`a l'identification des param\`etres}
\RRetitle{Reduced-order modeling for cardiac electrophysiology. Application to parameter identification}
\titlehead{Reduced-order modeling in cardiac electrophysiology}
\RRresume{
Un mod\`ele r\'eduit bas\'e sur la d\'ecomposition orthogonale aux valeurs propres (POD) est propos\'e pour les \'equations bidomaines de l'\'electro\-phy\-siologie cardiaque. Sa pr\'ecision est \'evalu\'ee \`a travers des \'electrocardiogrammes dans diverses configurations (mod\`eles d'infarctus du myocarde, si\-mu\-lations de longue dur\'ee, etc.). On montre en particulier que la courbe de restitution peut \^etre approch\'ee de mani\`ere efficace par cette approche. Le mo\-d\`e\-le r\'eduit est ensuite utilis\'e dans un probl\`eme inverse r\'esolu par un algorithme \'evolutionnaire. On pr\'esente des r\'esultats  d'identification de param\`etres ioniques et de localisation d'une zone infarcie \`a partir d'un ECG synth\'etique.}

\RRabstract{
A reduced-order model based on Proper Orthogonal Decomposition (POD) is proposed for the bidomain equations of cardiac electrophysiology. Its accuracy is assessed through electrocardiograms in various configurations, including myocardium infarctions and long-time simulations. We show in particular that a restitution curve can efficiently be approximated by this approach. The reduced-order model is then used in an inverse problem solved by an evolutionary algorithm. Some attempts are presented to identify ionic parameters and infarction locations from synthetic ECGs.}
\RRmotcle{\'Electrophysiologie cardiaque, mod\`eles r\'eduits, prob\`eme inverse, POD}
\RRkeyword{Cardiac electrophysiology, reduced-order model, inverse problem, POD}
\RRprojets{REO}
\RCParis 

\begin{document}
\RRNo{7811}

\makeRR   

\newtheorem{Rem}{Remarque}[section] 
\newcommand{\dep}[2]{\frac {\partial #1}{\partial #2}}
\newcommand{\Int}[2]{{\displaystyle \int_{#1}^{#2}}}
\newcommand{\Cm}{ C_{\mathrm {m}}   }
\newcommand{\OmHi}{\Omega_{\mathrm H,\mathrm i}}
\newcommand{\OmHe}{\Omega_{\mathrm H,\mathrm e}}
\newcommand{\OmH}{\Omega_{\mathrm H}}
\newcommand{\OmT}{\Omega_{\mathrm T}}
\newcommand{\OM}{\Omega}
\renewcommand{\div}{\textrm{div}\,}
\newcommand{\n}{\boldsymbol{n}}
\newcommand{\nT}{\boldsymbol{n}_{\mathrm T}}
\newcommand{\tr}{\mathrm{tr}}
\newcommand{\Jac}{\mathrm{Jac}}
\newcommand{\f}{\forall\,}
\renewcommand{\tr}{\mathrm{tr}}
\renewcommand{\l}{\lambda}
\newcommand{\e}{\epsilon}
\newcommand{\bb}[1]  {{\boldsymbol{#1} }}
\newcommand{\diverg}{ \operatorname{div} }
\newcommand{\grad}   { \bb{\nabla}    }
\newcommand{\sigmai}   { \bb{\sigma}_{\mathrm i}    }
\newcommand{\sigmab}   { \bb{\sigma}    }
\newcommand{\sigmae}   { \bb{\sigma}_{\mathrm e}    }
\newcommand{\sigmam}   { \bb{\sigma}_{\mathrm m}    }
\newcommand{\sigmaT}   { \bb{\sigma}_{\mathrm T}    }
\newcommand{\jt}   { \bb{j}    }
\newcommand{\ji}   { \bb{j}_{\mathrm i}    }
\newcommand{\je}   { \bb{j}_{\mathrm e}    }
\newcommand{\phii}   { \phi_{\mathrm i}    }
\newcommand{\phie}   { \phi_{\mathrm e}    }
\newcommand{\phiT}   { \phi_{\mathrm T}    }
\newcommand{\Gext}   { \Gamma_{\mathrm{ext}}}
\newcommand{\Gepi}   { \Gamma_{\mathrm{epi}}}
\newcommand{\Gendo}   {\Gamma_{\mathrm{endo}}}
\newcommand{\id}   { \bb{I}    }
\newcommand{\R}{\mathbb R}
\newcommand{\N}{\mathbb N}
\newcommand{\lapl}{\Delta}
\newcommand{\Iion}{ {I_{\mathrm {ion}}}   }
\def\eqd{\stackrel{\mathrm{def}}{=}}
\newcommand{\vm}   {v_{\mathrm m}}

\newcommand{\Frac}[2]{{\displaystyle \frac{#1}{#2}}}
\newcommand{\Der}[2]{\Frac {\partial #1}{\partial #2}}
\newcommand{\ue}{ u_{\mathrm e}   }
\newcommand{\Iapp}{ {I_{\mathrm {app}}}   }
\newcommand{\ui}{ u_{\mathrm i}   }
\newcommand{\Vm}{ V_{\mathrm m}   }
\newcommand{\Am}{ A_{\mathrm m}   }
\newcommand{\uj}{ u_{\mathrm j}   }
\newcommand{\uT}{ u_{\mathrm T}   }
\newcommand{\x}{\mathbf{x}}
\newcommand{\y}{\mathbf{y}}

\section{Introduction}

The bidomain equations used to model the cardiac electrophysiology is known to be very demanding from a computation viewpoint \cite{colli-pavarino-04,scacchi-pavarino-08}. To reduce its complexity, reduced-order models based on physical arguments can be considered. For example, assuming that the anisotropy ratios of the electrical conductivity tensors are the same in the intra and extracellular media, a simplified model called ``monodomain'' can be derived \cite{clements-nenonen-li-horacek-04,collifranzone-pavarino-taccardi-05}. Another example of physical simplification is given by the Eikonal model that essentially describes the propagation of the activation front \cite{keener-91,tomlinson-hunter-pullan-02}. Here we follow another route: keeping the physics of the model, we reduce the complexity of its discretization by using a Galerkin basis adapted to the kind of solutions that are looked for. This basis can be computed by different approaches. In this work, we choose the most popular one, which is known as the Proper Orthogonal Decomposition (POD) or Karhunen-Lo\`eve method.
This approach has been used in many fields of science and engineering, but to the best of our knowledge not for the equations of cardiac electrophysiology. The present study is a first step in this direction. We consider various configurations of interest in order to point out those where the reduced-order model seems promising and those where it has to be used with much care.  

After a brief presentation of the equations and the numerical methods used for the full and reduced-order methods in Section~\ref{sec:model-method}, the reduced-order model is tested in three situations in Section~\ref{sec:forward}: perturbation of some parameters, long-time simulation and infarction modeling. In all these cases, the results are assessed through the corresponding electrocardiogram (ECG). In Section~\ref{sec:inverse}, the reduced-order model is used with an evolutionary algorithm to identify some parameters of the problem and the location of infarcted regions from the electrocardiograms. The main conclusions of the study are then summarized in Section~\ref{sec:conclusion}.

\section{Models and methods}\label{sec:model-method}

\subsection{Models}\label{sec:model}

The electrical activity in the heart is modeled by the bidomain equations (see \cite{NK93,pennacchio-savare-colli-05,sachse,sundnes} e.g.). The extracellular potential $u_e$ and the transmembrane potential $\Vm$ are solutions in the heart domain $\OmH$ of the following system
\begin{equation}\label{eq:bidomain}
    \begin{aligned}
      A_{\mathrm m}\bigg(C_{\mathrm m} \Der{V_{\mathrm m}}{t} +I_{\mathrm{ion}}
      (V_{\mathrm m},w)\bigg)  -\diverg(\bb{\sigma}_{\mathrm i} \grad \Vm) 
- \diverg(\sigmai \grad \ue) =\Am  \Iapp\text{ in } (0,T)\times \OmH\\
      - \diverg( (\sigmai + \sigmae)  \grad \ue) - \diverg(\sigmai \grad \Vm)=0 \text{ in } (0,T)\times \OmH,
      \end{aligned}
\end{equation}
where $A_{\mathrm m}$ is the rate of membrane area per volume unit, $C_{\mathrm m}$ the membrane capacitance per area unit, $\sigmai$ and $\sigmae$ the intra- and extracellular conductivity tensors. \\
The term $\Iapp$ is a given source function, used in particular to initiate the activation, and $I_{\mathrm{ion}}(V_{\mathrm m},w)$ represents the ionic current across the membrane. The ionic variable $w$ is solution of the ordinary differential equation
\begin{equation}
 \frac{\partial w}{\partial t} +g(\Vm,w)=0 \text{ in } (0,T)\times \OmH.
\end{equation}

In this study, the dynamics of $w$ and $\Iion$ are described by the phenomenological two-variable model introduced by Mitchell and Schaeffer \cite{mitchell-schaeffer-03}:
\begin{equation}
  \label{eq:MS}
\begin{aligned}
  \Iion(\Vm,w) & = - \frac{ w}{\tau_{\mathrm{in}}}
  \frac{(\Vm-V_{\textrm{min}})^2(V_{\textrm{max}}-\Vm)}{V_{\textrm{max}}-V_{\textrm{min}}}
  + \frac{ 1} {\tau_{\mathrm{out}}} \frac{\Vm-V_{\textrm{min}}}{V_{\textrm{max}}-V_{\textrm{min}}},\\[2mm] 
   g(\Vm,w) & = 
  \left\{
    \begin{aligned}
      \Frac{w}{\tau_{\mathrm{open}}} -\Frac{1}{\tau_{\mathrm{open}}(V_{\textrm{max}}-V_{\textrm{min}})^2} &\quad \mbox{if}\quad \Vm < V_{\mathrm{gate}},\\ 
      \Frac{w}{\tau_{\mathrm{close}}} & \quad \mbox{if}\quad \Vm > V_{\mathrm{gate}},
    \end{aligned}
  \right.
\end{aligned}
\end{equation}
where $\tau_{\mathrm{in}}$, $\tau_{\mathrm{out}}$, $\tau_{\mathrm{open}}$,
$\tau_{\mathrm{close}}$, $V_{\mathrm{gate}}$  are given parameters and $V_{\mathrm{min}}$, $V_{\mathrm{max}}$ are scaling constants (typically -80 and 20 mV, respectively). Although this model describes in a very simplified way the ionic exchanges, it allows to recover realistic action potentials at the cell scale. Moreover, each parameter has a physiological meaning: the time constants  $\tau_{\mathrm{in}}$, $\tau_{\mathrm{out}}$
are respectively related to the length of the depolarization and repolarization (final stage) phases, $\tau_{\mathrm{open}}$ and 
$\tau_{\mathrm{close}}$ are the characteristic times of gate opening and closing respectively and $V_{\mathrm{gate}}$ corresponds to the change-over voltage.

On the boundary of the heart domain, it is generally admitted that the intracellular medium is isolated. So we have 
\begin{equation}\label{eq:insulating-intra}
\sigmai \grad \Vm \cdot \n +  \sigmai \grad \ue\cdot \n = 0\quad \text{on}\quad \partial \OmH.
\end{equation}
For the second boundary condition, two kinds of coupling between the heart and
the external medium -- that will be called the torso -- are usually
considered. The most complex one is to assume that the potentials and currents are
continuous between the extracellular and the torso, which corresponds
to a strong coupling. In the present study, the effect of the torso on
the heart is neglected and it is assumed that the extracellular medium is isolated (see  \cite{colli-pavarino-04,collifranzone-pavarino-taccardi-05}). This
corresponds to a weak coupling described by the following conditions: 
\begin{equation}\label{eq:couplage-faible}
  \left\{
    \begin{aligned}
     \sigmae\grad \ue\cdot \n &=0,&\quad\mbox{on}\quad \partial \OmH,\\
      \ue&=\uT,&\quad\mbox{on}\quad\partial \OmH.
    \end{aligned}
  \right.
\end{equation}

In the torso domain, denoted by $\OmT$, the electrical potential $\uT$ is solution of the generalized Laplace equation:
\begin{equation}\label{eq:thorax}
  \diverg(\sigmaT\grad \uT)=0,\quad \text{in}\quad \OmT.
\end{equation}
where $\sigmaT$ stands for the conductivity of the different tissues (lungs, bones, \emph{etc.}). 
Assuming that the body surface  $\Gamma_{\mathrm {ext}}$ is isolated, we have
\begin{equation}\label{eq:torso-gext}
\sigmaT\grad \uT\cdot \n=0,\quad\mbox{on}\quad\Gamma_{\mathrm {ext}}.
\end{equation}
This model allows to reproduce the main recording of the electrical cardiac activity made by medical doctor, the electrocardiogram (ECG).  It is a set of 12-graphs of voltages vs. time obtained from electrodes located on standard positions of the surface of the body. 

In \cite{boulakia-cazeau-fernandez-gerbeau-zemzemi-10}, under appropriate modeling assumptions (anisotropy, cell heterogeneity, simplified model of His bundle), realistic ECGs are provided in the healthy case and for two pathologies: left bundle branch block and right bundle branch block. The ECGs obtained satisfy the usual criteria used by medical doctors to detect these pathologies. We refer to  \cite{boulakia-cazeau-fernandez-gerbeau-zemzemi-10} for a description of the modeling choices made in our study. 

The heart-torso uncoupling conditions (\ref{eq:couplage-faible}) induces a significant reduction of the computational cost for the numerical resolution since the heart problem can be solved independently of the torso problem (in our case, the weak coupling is about 60 times faster than the strong coupling). 
As noticed in \cite{boulakia-cazeau-fernandez-gerbeau-zemzemi-10}, making this hypothesis increases the amplitude of the ECG but the shape of the ECG is only slightly modified. Thus, this simplification will be used in the present study.  \\

\subsection{Numerical methods for the full-order model}\label{sec:numresol}

To solve numerically the model given by (\ref{eq:bidomain})-(\ref{eq:torso-gext}), we follow the procedure explained in \cite{boulakia-cazeau-fernandez-gerbeau-zemzemi-10} with the heart-torso uncoupling assumption. The spatial discretization  is based on the finite element method and the time discretization is performed using the second order BDF implicit scheme with an explicit treatment of the ionic current. \\

Let us first introduce some notations to present the resolution algorithm. Let $V_h \subset H^1(\OmH)$ and $W_h \subset H^1(\OmT)$ be two finite dimensional subspaces of continuous piecewise affine functions. We define $W_{h,0}=\big\{�\psi_{\mathrm T} \in W_h / \psi_{\mathrm T}=0 \text{ on } \partial\OmH\big\}$. \\
We assume that the time interval $[0,T]$ is divided in $K$ subintervals $[t_k,t_{k+1}]$, $0\leq k\leq K-1$, with $t_k\eqd k\delta t$ where $\delta t \eqd T/K $ is the time step. We denote by $(\Vm^{k},\ue^k,w^{k})$ the approximated solution
obtained at time $t_k$. Then, $(\Vm^{k+1},\ue^{k+1},w^{k+1})$ is computed as follows: 

\begin{itemize}
\item First, a second order extrapolation of $\Vm$ at time $t_{k+1}$ and $w^{k+1}\in V_h$ are computed: \\ $\widetilde V_{\mathrm m}^{k+1}=2\Vm^{k}-\Vm^{k-1}$ and 
$$ \frac{1}{\delta t} \left(\frac{3}{2} w^{k+1} - 2 w^{k} +\frac{1}{2} w^{k-1} \right)+g(\widetilde V_{\mathrm m}^{k+1},w^{k+1})=0$$
	\item Then, $(\Vm^{k+1},\ue^{k+1})\in V_h\times V_h$, with $\displaystyle{\int_{\OmH}\ue^{k+1} = 0}$ is computed as the solution of:
\begin{equation}\label{heart-scheme}
  \left\{
    \begin{aligned}
     & \begin{split}
\Am \Int{\OmH}{} \frac{\Cm}{\delta t} \left(\frac{3}{2} \Vm^{k+1} - 2 \Vm^{k} +\frac{1}{2} \Vm^{k-1} \right)\phi +    
\Int{\OmH}{}        \sigmai  \grad (\Vm^{k+1 } + \ue^{k+1 }) \cdot \grad \phi \\
       =  \Am \Int{\OmH}{}  \left( \Iapp(t_{k+1})  - \Iion({\widetilde V_{\mathrm m}^{k+1}},w^{k+1}) \right)\phi,
\end{split}\\
 & \int_{\OmH} \sigmai \grad\Vm^{k+1}\cdot\grad\psi_{\mathrm e}+\int_{\OmH} (\sigmai +\sigmae) \grad \ue^{k+1}\cdot\grad\psi_{\mathrm e}  = 0,
\end{aligned}
\right.
\end{equation}
for all $(\phi,\psi_{\mathrm e})\in V_h \times V_h$, with $\displaystyle{\int_{\OmH}\psi_{\mathrm e} = 0}$.
\end{itemize}
For other strategies for the time discretization, we refer to \cite{austin-trew-pullan-06,colli-pavarino-04,ethier-bourgault-08,lines-grottum-tveito-03,sundnes-lines-tveito-01} for semi-implicit schemes or \cite{keener-bogar-98,sundnes-lines-tveito-05,vigmond-weber-et-al-08} for operator splitting approaches. 

Due to the second order approximation $\widetilde V_{\mathrm m}^{k+1}$ of $\Vm^{k+1}$, the proposed time scheme has the convenient property to be associated to a constant matrix. Indeed, since the ionic current is explicit,  (\ref{heart-scheme}) is equivalent to a linear system whose matrix in $\mathbb{R}^{2N, 2N}$ is given by
\begin{equation}\label{sys-matrix}
\begin{pmatrix}
\dfrac 32\dfrac{\Am\Cm}{\delta t} \bb{M}+\bb{K_1}& \bb{K_1}\\
{ } \\
\bb{K_1} & \bb{K_2}
\end{pmatrix}.
\end{equation}
where matrices $\bb{M}  \in \mathbb{R}^{N, N}$,  $\bb{K_1} \in \mathbb{R}^{N,N}$ and $\bb{K_2} \in \mathbb{R}^{N, N}$ are given by
$$
\bb{M}_{ij}= \Int{\OmH}{ } \phi_i\phi_j,\,(\bb{K_1})_{ij}= \Int{\OmH}{} \sigmai \grad  \phi_i\cdot\grad  \phi_j,\,(\bb{K_2})_{ij}=\Int{\OmH}{}  (\sigmai +\sigmae) \grad  \phi_i\cdot\grad  \phi_j
$$
where $(\phi_i)_{1 \leq i \leq N}$ is a finite element basis of $V_h$.

Once $u_e^{k+1}$ is known, the torso potential $\uT^{k+1} \in W_{h}$ is obtained by solving:
\begin{equation}
\left\{
\begin{aligned}
\int_{\OmT} \sigmaT \grad \uT^{k+1}\cdot\grad\psi_{\mathrm T} = 0,& \quad \forall \psi_{\mathrm T} \in W_{h,0}\\
\uT^{k+1}
      =\ue^{k+1},&\quad\mbox{on}\quad\partial \OmH.
\end{aligned}
\right.
\end{equation}\label{torso-scheme}
Since this last step depends linearly on $\ue^{k+1}$ and is weakly coupled to the rest of the problem, it can be solved very efficiently by precomputing a transfer matrix \cite{boulakia-cazeau-fernandez-gerbeau-zemzemi-10}. Thus, the main computational cost comes from the resolution of system (\ref{heart-scheme}) satisfied by $(\Vm^{k+1},\ue^{k+1})$.

\subsection{Reduced-order Modeling}\label{section:POD}

For completeness, we briefly recall the principle of the the Proper Orthogonal Decomposition (POD) method. We refer the reader interested by more details to \cite{kunisch-volkwein-01,rathinam-petzold-04} for example. POD is a method used to derive low-order models by projecting the system onto subspaces spanned by a basis of elements that contains the main features of the expected solution. To generate the POD basis associated with a precomputed solution $u=(\Vm,\ue)$ of the Galerkin problem (\ref{heart-scheme}), we make a first numerical simulation (or a set of simulations) and keep some snapshots $u(t_k), \, 1 \leq k \leq p$. Then a singular value decomposition (SVD) of the matrix $B=(u(t_1),\dots,u(t_p)) \in \mathbb{R}^{2N,p}$ is performed:
$$B=USV',$$ 
where $U \in \mathbb{R}^{2N,2N}$ and $V\in \mathbb{R}^{p,p} $
are orthogonal matrices, $S\in \mathbb{R}^{2N,p}$ is the matrix of the singular values ordered by decreasing order. The $N_{\mathrm{modes}}$ first POD basis functions $\{\Psi_i \}_{1 \leq i \leq N_{\mathrm{modes}}}$ are then given by the $N_{\mathrm{modes}}$ first columns of $U$ and the POD Galerkin problem is solved by looking for a solution of the type $$u=\displaystyle\sum_{i=1}^{N_{\mathrm{modes}}}\alpha_i(t)\Psi_i.$$
The $2N \times 2N$ sparse matrix  given by (\ref{sys-matrix})  is thus replaced by a full matrix of size $N_{\mathrm{modes}} \times N_{\mathrm{modes}}$ with the POD method. Since this matrix is constant in time, it is therefore projected on the POD basis and factorized (or even inversed) only once at the beginning of the computation. For the simulations presented in this paper, the order of magnitude of $N_{\mathrm{modes}}$ is typically one hundred and the reduced-order model resolution is about one order of magnitude faster than the full-order one. \\

Note that the model is reduced by changing the basis the solution is approximated on, but it still corresponds to the discretization of the original system~\eqref{heart-scheme}. In particular, it depends on exactly the same parameters as the full-order model. In general, the reduced-order model does very well at reproducing the solutions that have been used to generated the POD basis. But it is difficult to anticipate how it behaves when the parameters are modified. This issue is difficult and is still the topic of active researches, see for example~ \cite{amsallem-farhat-08}, \cite{chaturantabut-sorensen-10}, \cite{grepl-maday-nguyen-patera-07}. In the next sections, we will show situations where the reduced model usually works correctly and others where it does not. For the latter, we propose some strategies to improve the results. 

\section{Application to forward problems}\label{sec:forward}

In this section, the reduced-order model is used in different configurations and compared to the full-order one.
It is important to note that the accuracy is not assessed on the whole solution, but only on the ECG, which is considered as the ``output of interest'' of these simulations.

\subsection{POD when some model parameters vary}\label{sec:pod-ex}

A POD basis is computed for a given set of parameters. We consider in this section what happens when this basis is used to solve a reduced-order model with different parameters. Our investigation is restricted to $\tau_{\rm in}, \tau_{\rm close}$, $A_{\rm m}$ and $C_{\rm m}$ because it has been shown in \cite{boulakia-cazeau-fernandez-gerbeau-zemzemi-10} that they are the ones the ECG is the most sensitive to. All the simulations are performed on an idealized geometry of ventricles based on two ellipsoids (see \cite{sermesant-et-al-06}). The mesh contains 418465 tetrahedra, the time-step is fixed to 0.5 milliseconds and the parameters of the model (\ref{eq:bidomain})-(\ref{eq:MS}) are given in Table \ref{tab:cell-membr-parameters}.

\begin{table}[htbp]
\centering
\[
\begin{array}{c|c|c|c|c|c|c|c|c|c|c|c}
 A_{\mathrm {m}} \mathrm{(cm^{-1})} & C_{\mathrm {m}}\mathrm{(mF)} & \tau_{\mathrm {in}} & \tau_{\mathrm {out}} & \tau_{\mathrm{open}} & \tau_{\mathrm{close}}^{\mathrm{RV}} & \tau_{\mathrm{close}}^{\mathrm{endo}} & \tau_{\mathrm{close}}^{\mathrm{mcell}} & \tau_{\mathrm{close}}^{\mathrm{epi}} & V_\mathrm{gate} & V_\mathrm{min} & V_\mathrm{max} \\
\hline
200 & 10^{-3} &  16 & 360 & 100 & 120 & 130 & 140 & 90 & -67 & -80 & 20
\end{array}
\]
\caption{Cell membrane parameters.}
\label{tab:cell-membr-parameters}
\end{table}

\medskip

\begin{figure}[h!]
  \centering
 \hspace*{-0.1cm}\includegraphics[width=0.4\linewidth]{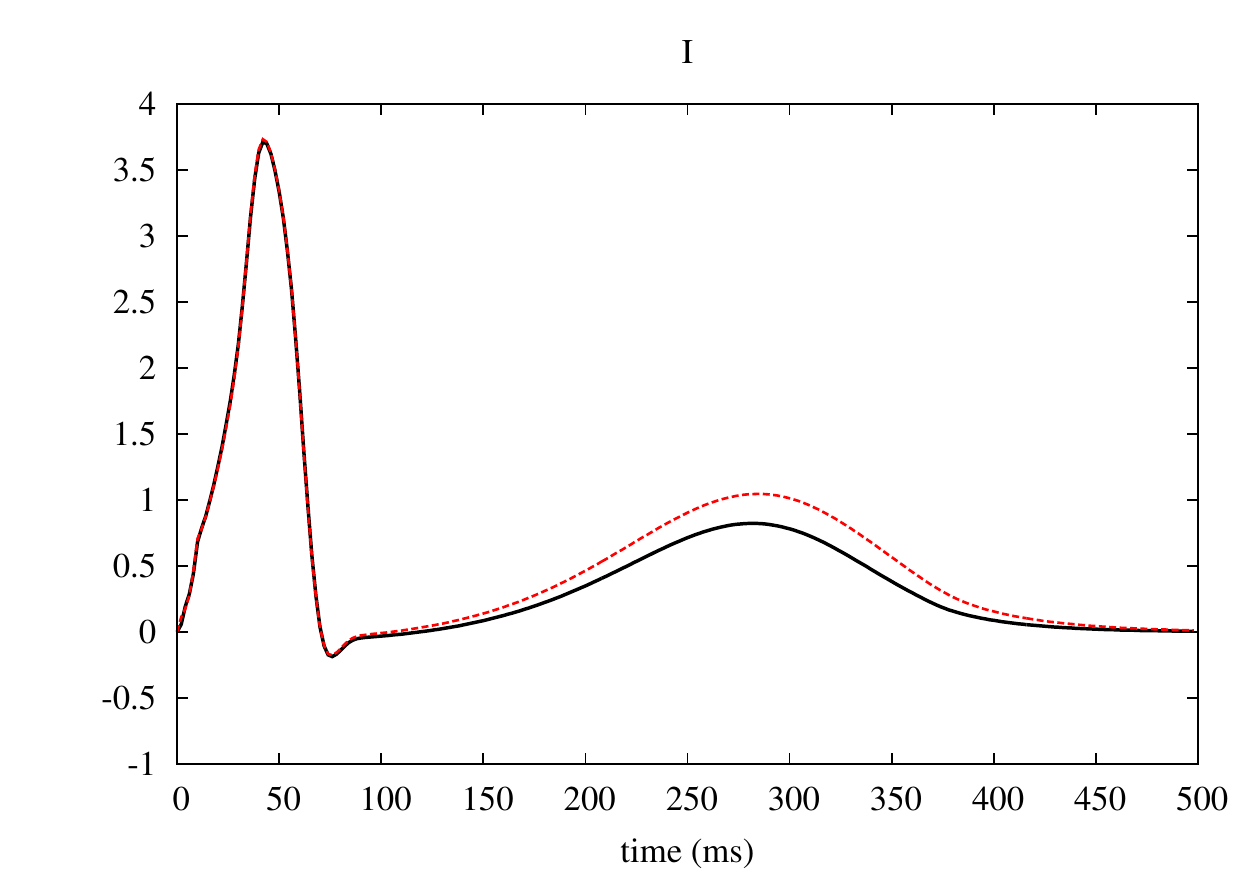}\hspace{0.2cm}
\includegraphics[width=0.4\linewidth]{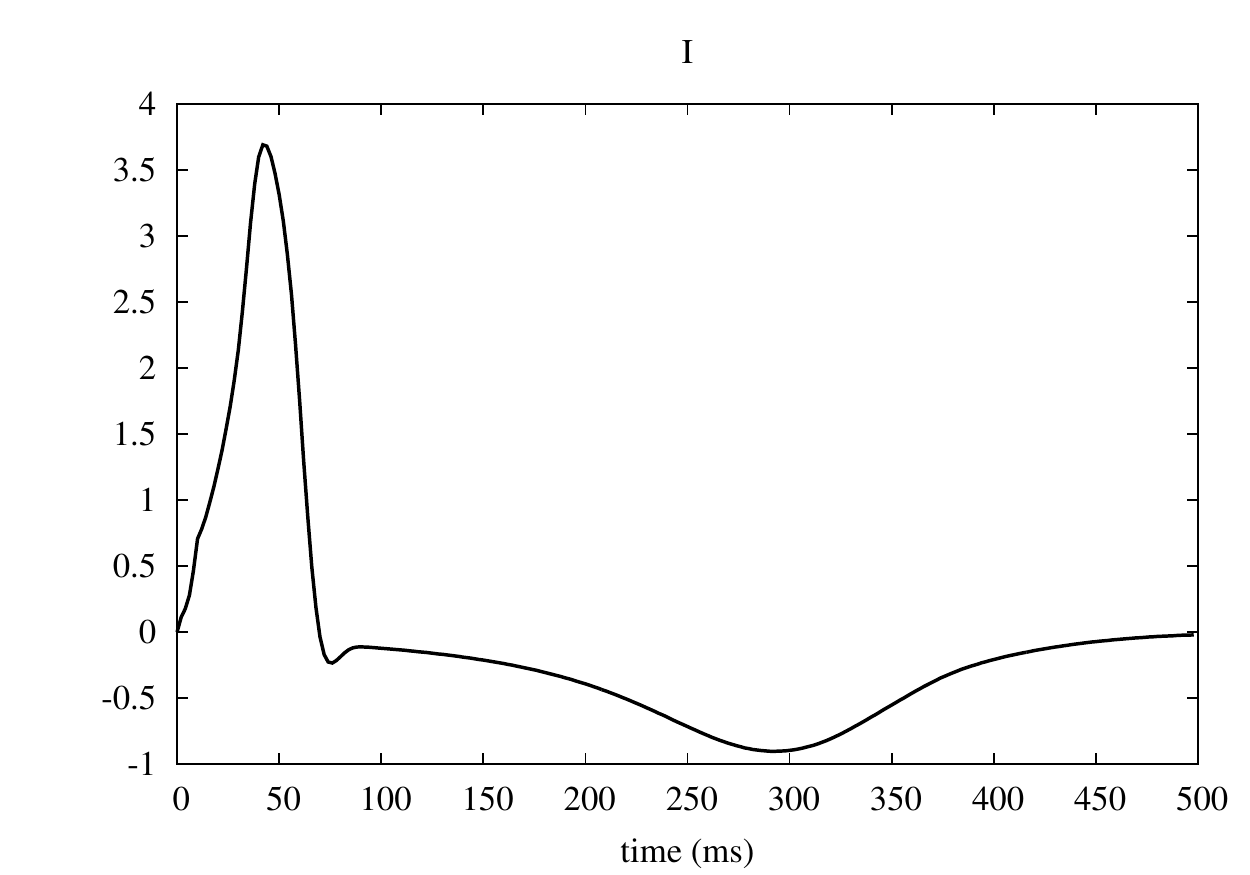}
 \caption{Left: First lead of the ECGs with $(\tau_{\mathrm {close}}^{\mathrm {epi}},\tau_{\mathrm
     {close}}^{\mathrm{RV}})=(80,130)$. Comparison between the full model (dotted red line)
 and the reduced model (black line) with the POD basis generated with $(\tau_{\mathrm {close}}^{\mathrm {epi}},
\tau_{\mathrm {close}}^{\mathrm {RV}}) =(100,100)$. 
     Right: First lead of the ECG computed with the full-order model for $(\tau_{\mathrm {close}}^{\mathrm {epi}},
\tau_{\mathrm {close}}^{\mathrm {RV}}) =(100,100)$.}
  \label{fig:ECGtau-close}
\end{figure}

Let us first consider a perturbation of $\tau_{\mathrm {close}}$, which mainly affects the repolarization phase \emph{i.e.}  the T-wave in the ECG.  In our model, the heart is divided in four regions where $\tau_{\mathrm {close}}$ takes different constant values (see \cite{boulakia-cazeau-fernandez-gerbeau-zemzemi-10}). We consider the value in the epicardium of the left ventricle $\tau_{\mathrm {close}}^{\mathrm {epi}}$ and in the right ventricle
$\tau_{\mathrm {close}}^{\mathrm {RV}}$. A POD basis of 80 vectors
is first constructed with $(\tau_{\mathrm {close}}^{\mathrm {epi}},\tau_{\mathrm {close}}^{\mathrm
  {RV}}) =(100,100)$.  
Then the corresponding reduced-order model is solved with $(\tau_{\mathrm {close}}^{\mathrm {epi}}, \tau_{\mathrm {close}}^{\mathrm {RV}}) = (80,130)$. Figure \ref{fig:ECGtau-close} (left) shows that the full and reduced models are in good agreement. We only refer to the first lead of the ECG for convenience, but the same trend is observed on the other leads.  It is interesting to note that the experiment used to generate the POD basis (with $(\tau_{\mathrm {close}}^{\mathrm {epi}}, \tau_{\mathrm {close}}^{\mathrm {RV}}) =(100,100)$) has a negative T-wave (Figure \ref{fig:ECGtau-close}, right). It is therefore not trivial that the reduced model is able to provide the correct (positive) T-wave.  If, instead of considering the ECG only, we compare the extracellular potential $u_e$ of the full and reduced models, the relative difference is $2.3 \times 10^{-2}$, in Euclidean norm in $\OmH\times[0,T]$, with $T=500 ms$. The full 3D solutions are therefore in good agreement too. \\

\begin{figure}[h!]
  \centering
\includegraphics[width=0.4\linewidth]{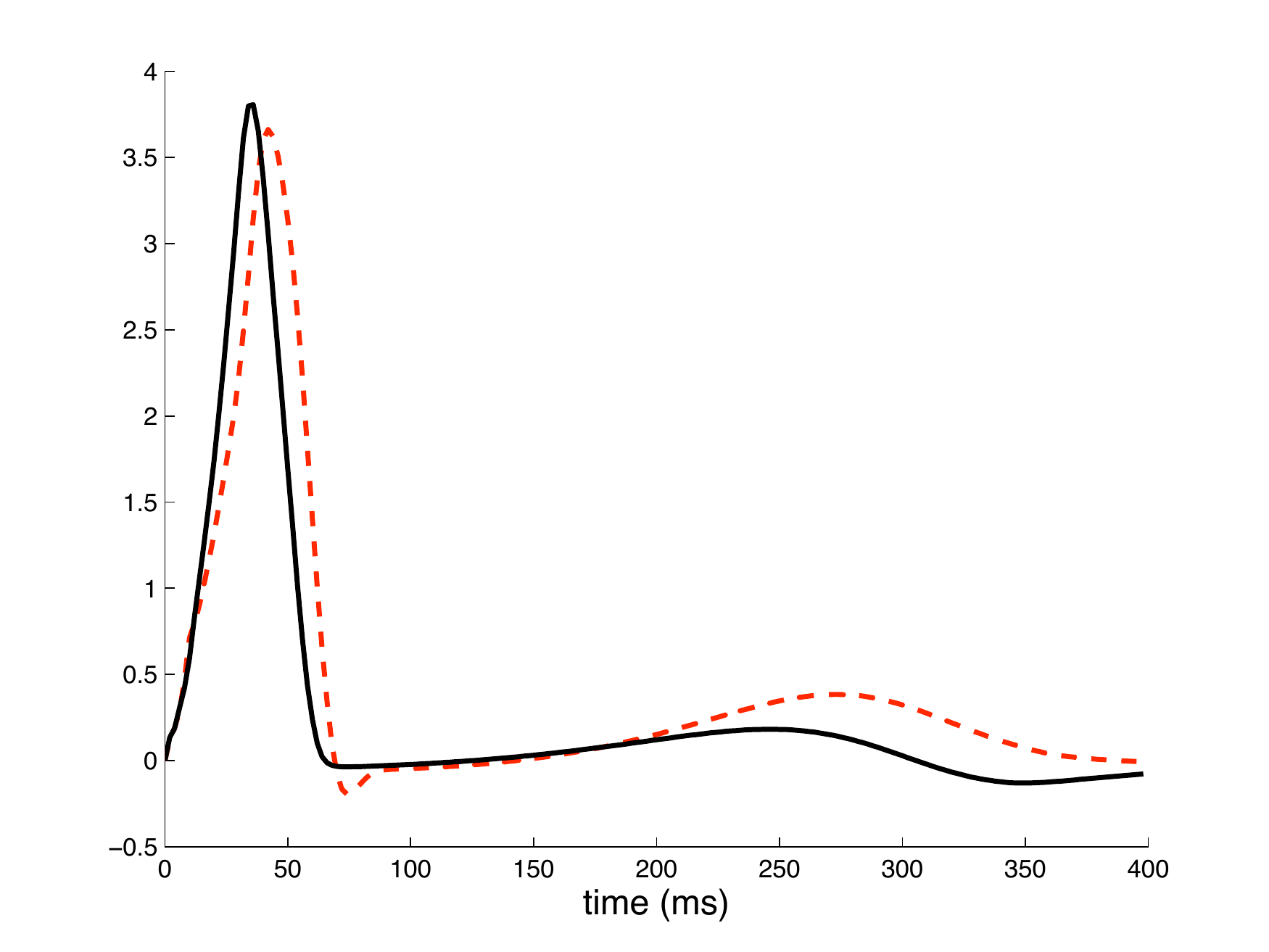}\hspace{-0.2cm}
\includegraphics[width=0.4\linewidth]{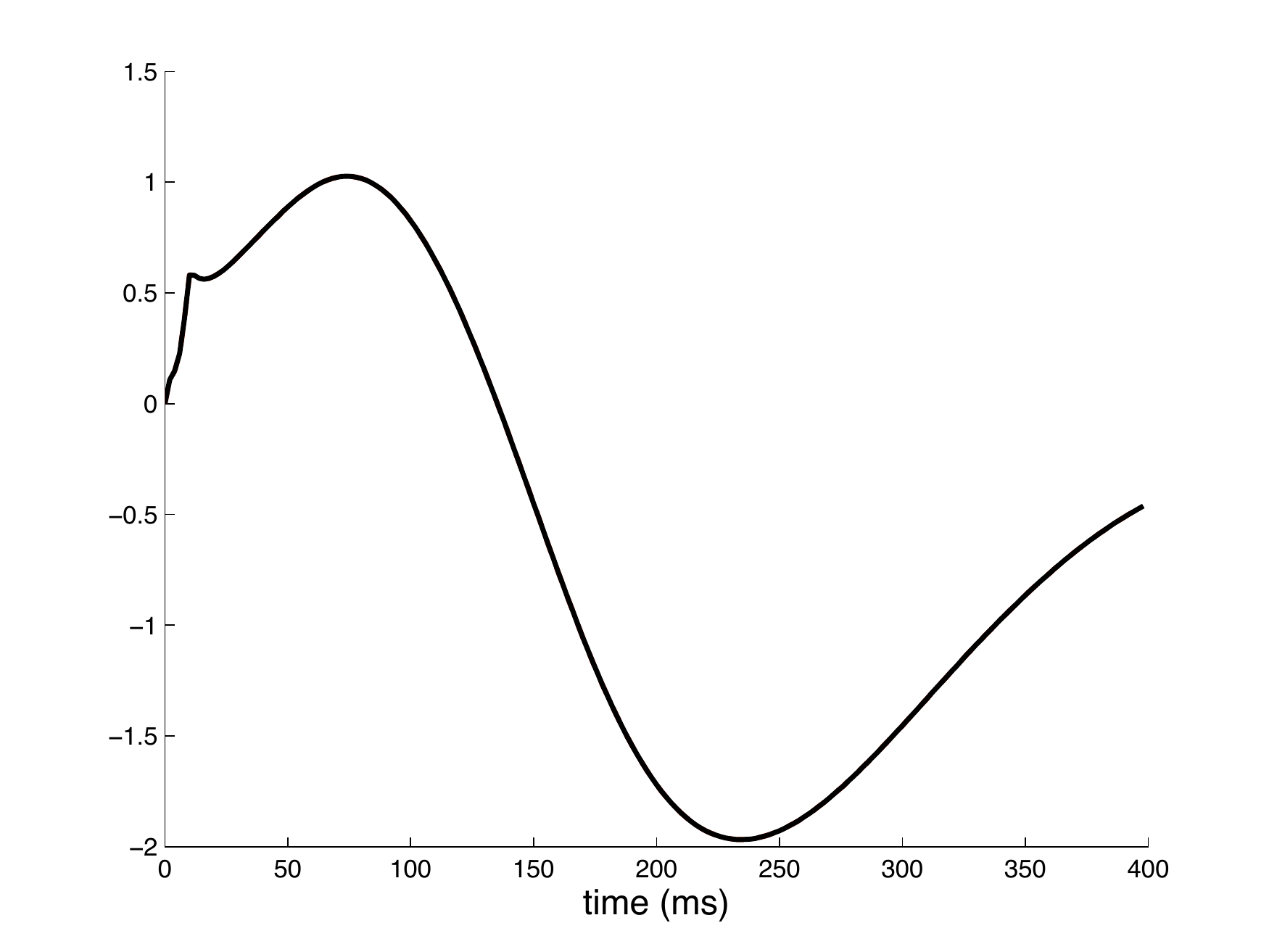}
 \caption{Left: First lead of the ECGs with $(\tau_{\mathrm{in}},\,C_{\mathrm{m}}, \,A_{\mathrm{m}},\,\tau_{\mathrm {close}}^{\mathrm{RV}})=(0.8,10^{-3}, 200, 120)$. 
Comparison between the full model (dotted red line)
 and the reduced model (black line) with a POD basis generated with $(\tau_{\mathrm{in}},\,C_{\mathrm{m}}, \,A_{\mathrm{m}},\,\tau_{\mathrm {close}}^{\mathrm{RV}})=(1.5,2\times 10^{-3},100,50)$. Right : First lead of the ECG computed with the full-order model for 
 $(\tau_{\mathrm{in}},\,C_{\mathrm{m}}, \,A_{\mathrm{m}},\,\tau_{\mathrm {close}}^{\mathrm{RV}})=(1.5,2\times 10^{-3},100,50)$}
  \label{fig:ECG4param}
\end{figure}  

Consider now perturbations of  $\tau_{\mathrm{in}}$, $C_{\mathrm{m}}$, $A_{\mathrm{m}}$ and $\tau_{\mathrm {close}}^{\mathrm  {rv}}$. Simulations are run with $\tau_{\mathrm{in}}=0.8,
\,C_{\mathrm{m}}=10^{-3}, \,A_{\mathrm{m}}=200, \,\tau_{\mathrm {close}}^{\mathrm
  {rv}}=120$. We compare the ECGs obtained with the full-order model and with
the reduced-order model corresponding to a POD basis generated with $\tau_{\mathrm{in}}=1.5,
\,C_{\mathrm{m}}=2 \times 10^{-3}, \,A_{\mathrm{m}}=100, \,\tau_{\mathrm {close}}^{\mathrm
  {rv}}=50$. In spite of some differences (slight temporal shift or difference of amplitudes, see Figure \ref{fig:ECG4param}, left), the results reasonably match.
Again, this matching is not trivial since the ECG corresponding to the values used to generated the POD basis  
 looks totally different (see Figure \ref{fig:ECG4param}, right). 

\medskip

Although it is not possible to make a general statement, these results are representative of several experiments that have been performed and not reported here:  the reduced-order model usually gives an acceptable ECG when the parameters $\tau_{\rm in}, \tau_{\rm close}$, $A_{\rm m}$ and $C_{\rm m}$ vary within a reasonable range. When these parameters vary too much, it is recommended to generate new POD bases, closest to the region of interest. This will be done in Section~\ref{sec:inverse} to address inverse problems.

\subsection{POD for long-time simulations and restitution curves}\label{sec:long-time}

The \emph{restitution curve} represents the dependence of the action potential duration (APD) on the preceding diastolic interval (DI). This curve is known to be important in the understanding of some arrhythmias~\cite{qu-xie-weiss-10,mitchell-schaeffer-03}. In addition, it has been shown in~\cite{manriquez-zhang-sorine-08} that the restitution curve can be used to identify some parameters in (\ref{eq:MS}). The simulation of the restitution curve is extremely challenging for 3D models since it requires several dozen of heart beats, whereas the simulation of one heart beat is already very demanding. We propose to investigate reduced-order modeling to simulate a series of heart beat with variable durations.  

\begin{figure}[!ht]
 \centering
  \includegraphics[width=0.98 \textwidth]{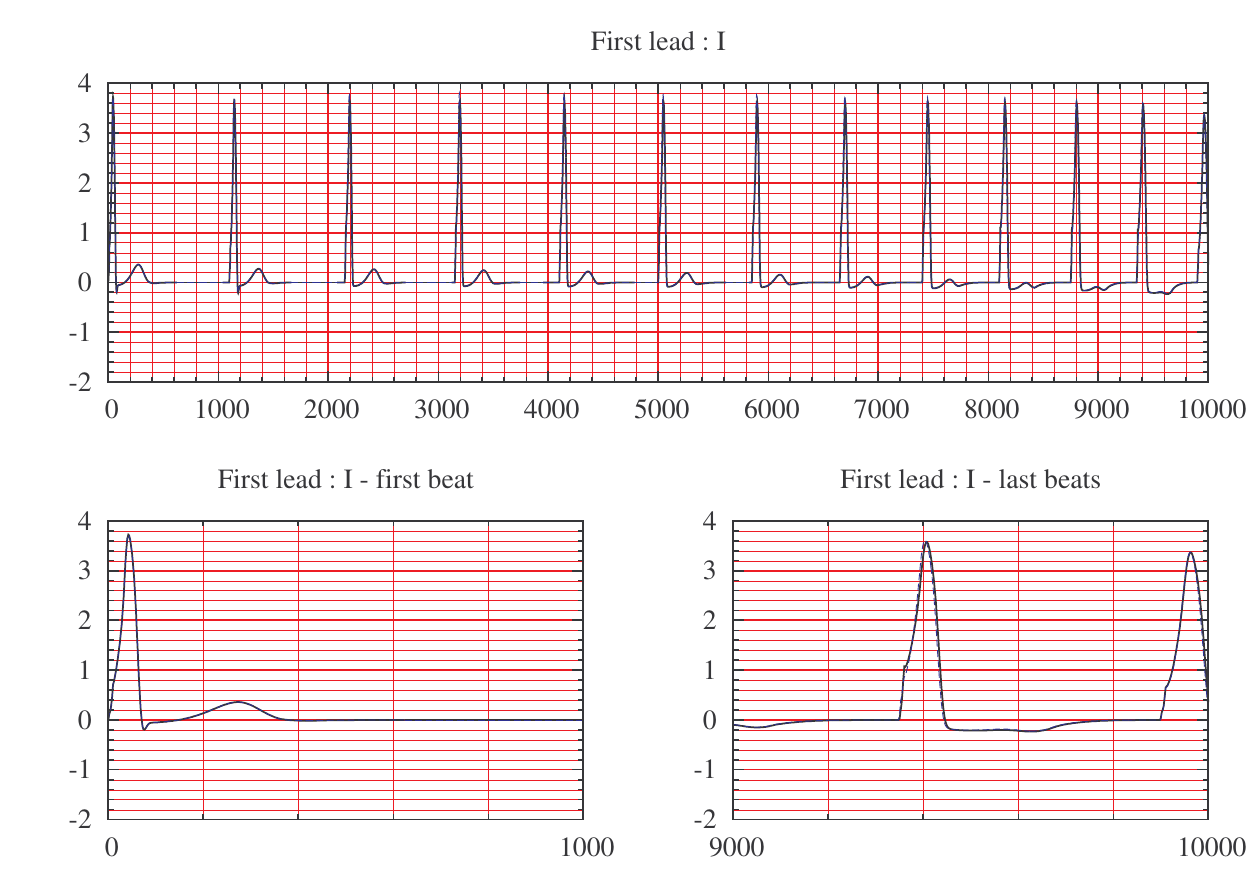}
  \caption{First lead of the ECG over a long period ($10$ seconds) for an increasing heart rate ($55$ beats per minute to $110$ beats per minute). On the left-hand side, a zoom of the first second of the simulation, on the right-hand side a zoom of the last second of the simulation. Solid line indicates the complete model solution, dotted line the reduced model one (they are almost superimposed).}
\label{fig:ecg_long} 
\end{figure}

\begin{figure}[!ht]
 \centering
  \includegraphics[width=0.98 \textwidth]{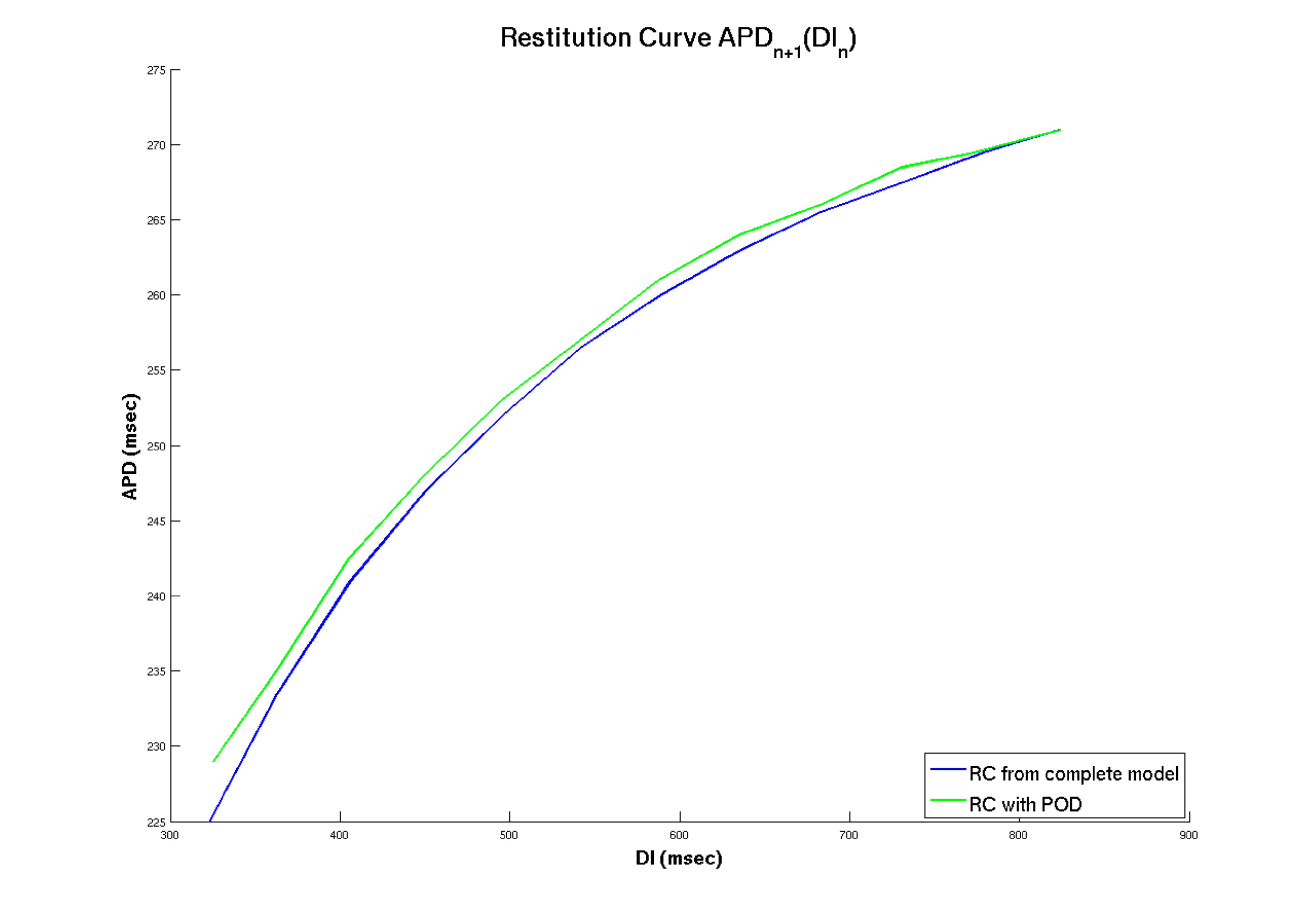}
  \caption{Restitution curve obtained solving a 12 second simulation with the complete model (blue line)  and the POD (green line).}
\label{fig:restitution_curve} 
\end{figure}

We consider a sequence of 12 heart beats, with a heart rate going from about 55 bpm to about 110 bpm (at each beat, the period is decreased of 50 ms). The simulation is first run with the full-order model for 400 milliseconds (less than one heart beat) with a timestep of 0.5 milliseconds. A snapshot is taken every 4 time steps. Then a POD basis of 100 modes is constructed from these snapshots and this basis is used for the long-time simulation (10 seconds). For comparison purpose, the full-order model was also run for the long-time simulation. Figure \ref{fig:ecg_long} shows the first lead of the ECG obtained with the full and reduced-order models: the two solutions are almost superimposed during the whole simulation. To take a closer look, we plot in Figure~\ref{fig:restitution_curve} the corresponding restitution curves obtained by recording the APD and DI in a cell localized at the epicardium of the right ventricle. The two curves are in very good agreement. The reduced-order model therefore seems to be able to correctly handle a heart rate acceleration, even when the POD basis has been generated from only one heart beat.  Of course, if the shape of the propagation front was significantly altered because of the heart beat increase, the results would probably deteriorate. But our results show that within a significant range of heart rates the reduced-order model is doing very well.

\subsection{Full and reduced-order simulation of an infarction}\label{sec:PODinfarct}

In \cite{boulakia-cazeau-fernandez-gerbeau-zemzemi-10}, the ECG model has been used to simulate healthy cases and bundle branch blocks. We propose in the present work to see the effect on the ECG of myocardial infarctions. An infarcted region cannot be activated anymore because it has been damaged by an interruption of blood supply. In our simplified framework, this defect is modeled by dividing by 10 the parameter $\tau_{\mathrm{out}}$ in the infarcted region. Doing so, the region keeps unactivated even when it is stimulated by the action potential. With more sophisticated ionic models, we could for example modify the behavior of the extracellular potassium~\cite{potse-coronel-leblanc-07}.

According to electrocardiology books \cite{dubin-03,PCFLV-07}, the main consequence of an infarction on a real ECG is an elevation or a depression of the ST segment in different leads. The magnitude of this elevation or depression and the leads where theses changes are visible depend on the position of infarction. Particularly, the main features we should find are:
\begin{itemize}
 \item in the case of a posterior infarction: a depression in the ST segment in the $V_1$ and $V_2$ leads;
 \item in the other cases: an elevation in the ST segment with an inverted T wave;
 \item in the case of an anterior infarction we should look at $V_1,V_2$ or $V_3$, in the case of a lateral infarction at $\text{I}$ or $aVL$ and for an inferior infarction at $\text{II}, \text{III}$ or $aVF$.
\end{itemize}

\begin{figure}[!htp]
 \centering
  \includegraphics[width=0.9 \textwidth]{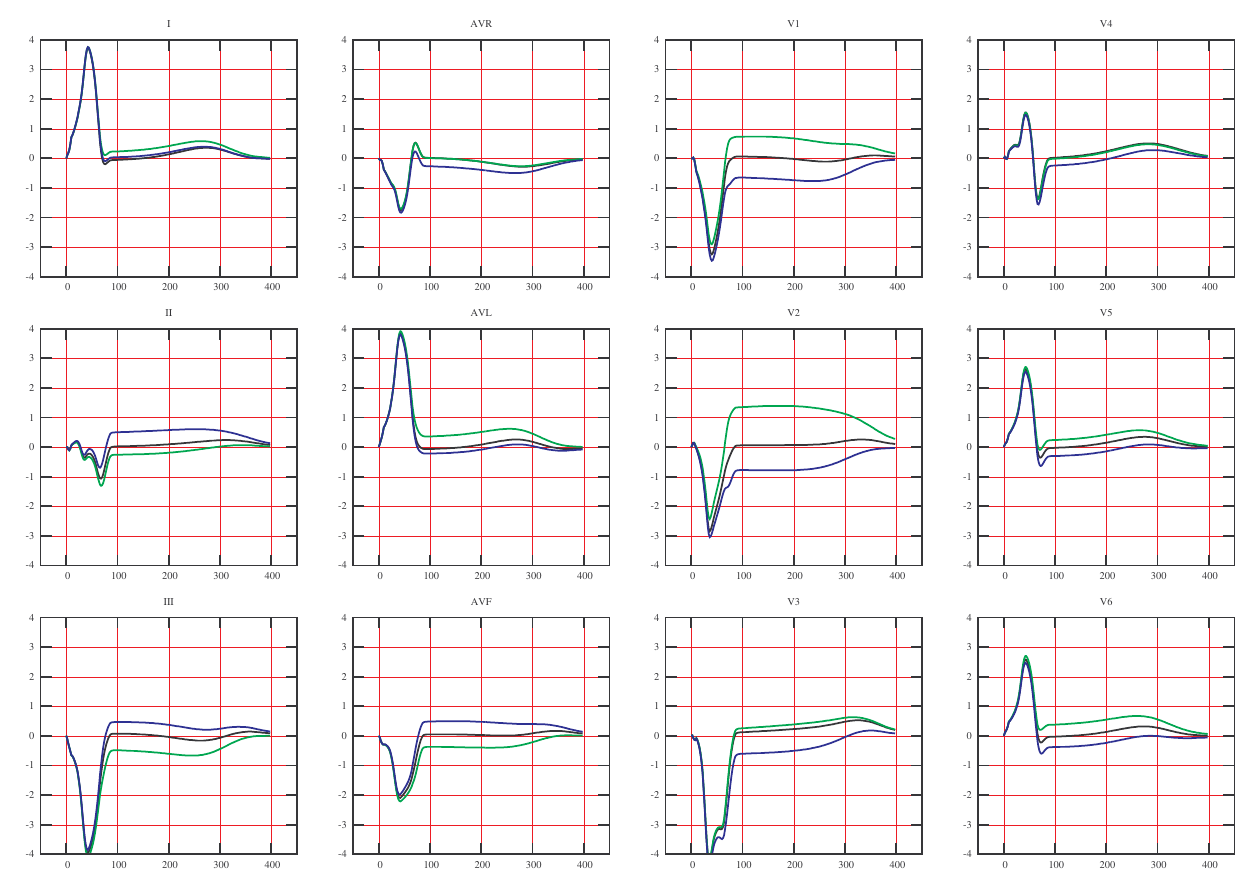}
  \caption{Simulated ECG: black line represents the healthy case, green line the transmural myocardial anterior infarction and blue one the posterior infarction.}
\label{fig:ecg4infarctions1} 
\end{figure}

\begin{figure}[!htp]
 \centering
  \includegraphics[width=0.9 \textwidth]{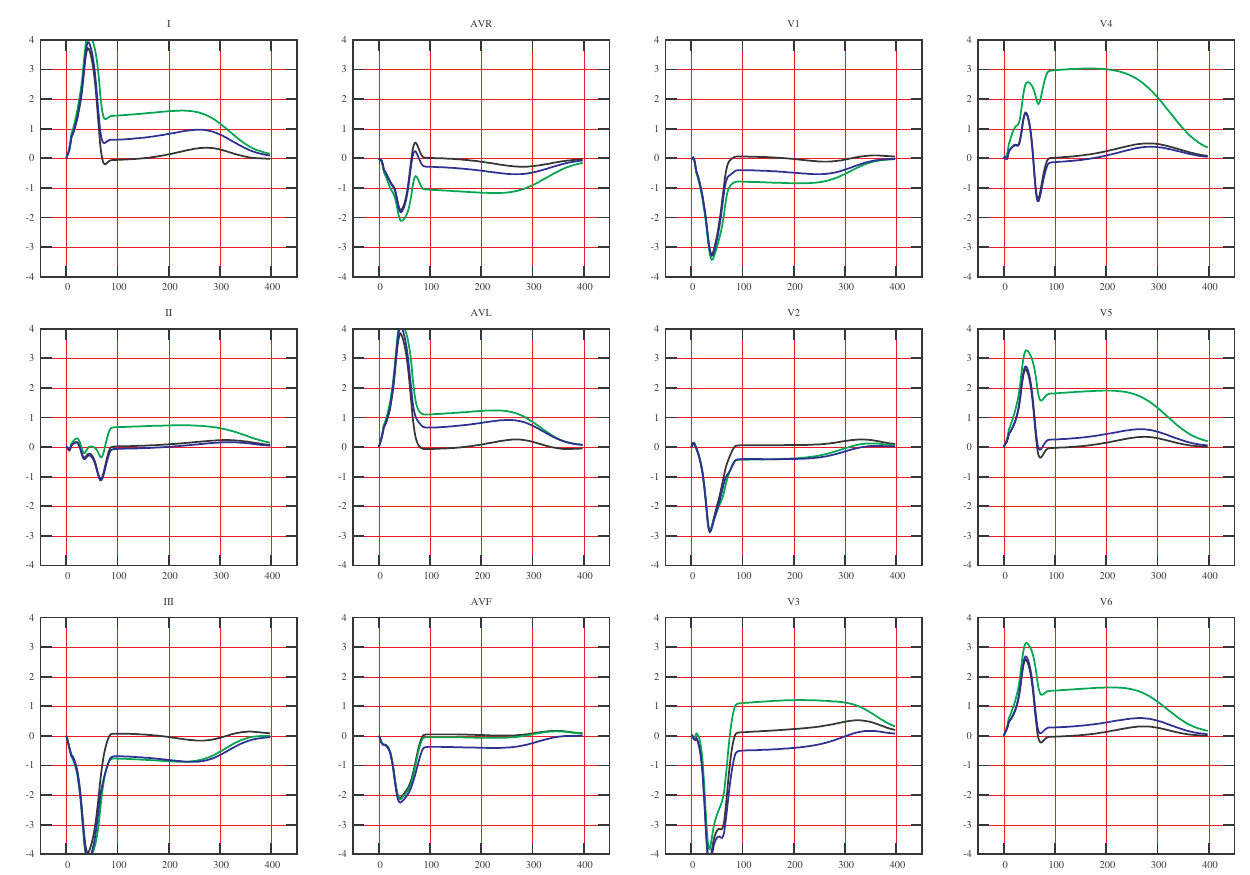}
  \caption{Simulated ECG: black line represents the healthy case, green line the transmural myocardial inferior infarction and blue one the lateral infarction.}
\label{fig:ecg4infarctions2} 
\end{figure}

Figures \ref{fig:ecg4infarctions1} and \ref{fig:ecg4infarctions2} show the simulated ECG for the main kinds of infarction: in Figure \ref{fig:ecg4infarctions1} the anterior and the posterior ones, and in Figure \ref{fig:ecg4infarctions2} the lateral and the inferior ones. The result are not very good for the inferior infarction: as expected we see an ST elevation in $\text{II}$ but we also observe an important ST elevation in the last three leads, a depression in $\text{III}$ and no sign on $aVF$ which is not expected. The difficulty of simulating the inferior infarction is probably due to the fact that this zone is very close to the initial activation region. It also seems that the QRS complex are not exactly as they should be. Nevertheless, we find as expected a depression (resp. elevation) of the ST segment in $V_1$, $V_2$ and $V_3$ leads in presence of a posterior (resp. anterior) infarction. We also find a ST elevation in~$\text{I}$ and~$aVL$ in presence of a lateral infarction. In conclusion, the results obtained with the full-order model for the ST segments are very satisfactory for the anterior, posterior and lateral infarctions.

\medskip

We now propose to investigate the same situations with the reduced-order model. Contrary to the conclusions of the previous two sections, we will see that a too naive use of POD gives in this case very poor results. 

\begin{figure}[htbp]
\begin{minipage}{0.42\textwidth}
\centering
\includegraphics[viewport = 0 0 849 640 , width=0.9\textwidth]{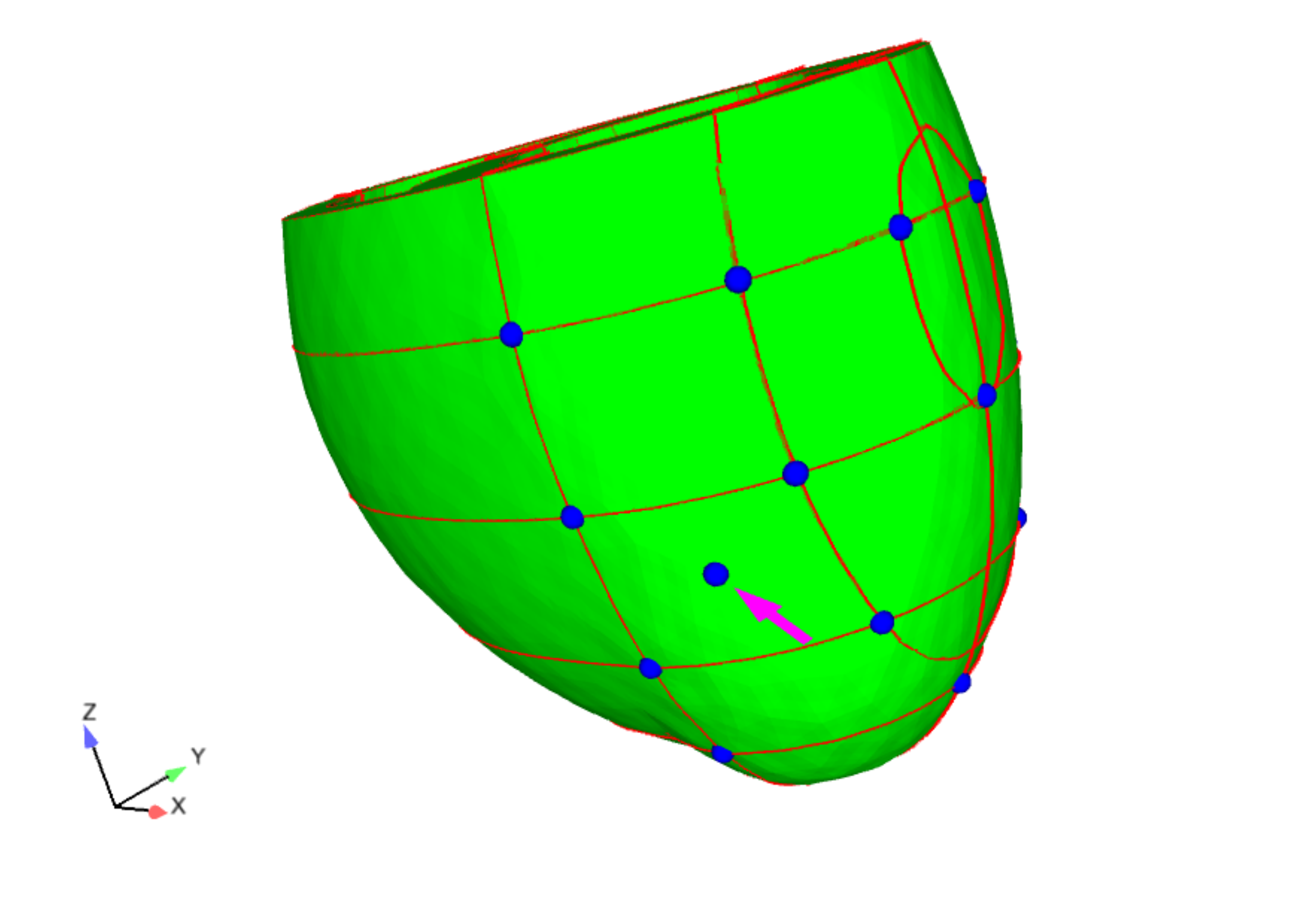}
\caption{The mesh and the 18 points used to build the infarction POD basis. The point out of the mesh lines, indicated by the pink arrow, is the point $P$ considered in the example.}
\label{fig:inf_points}
\end{minipage} \hspace{0.05 cm} \ \begin{minipage}{0.55\textwidth}
                  \centering
                  \includegraphics[width=1\linewidth]{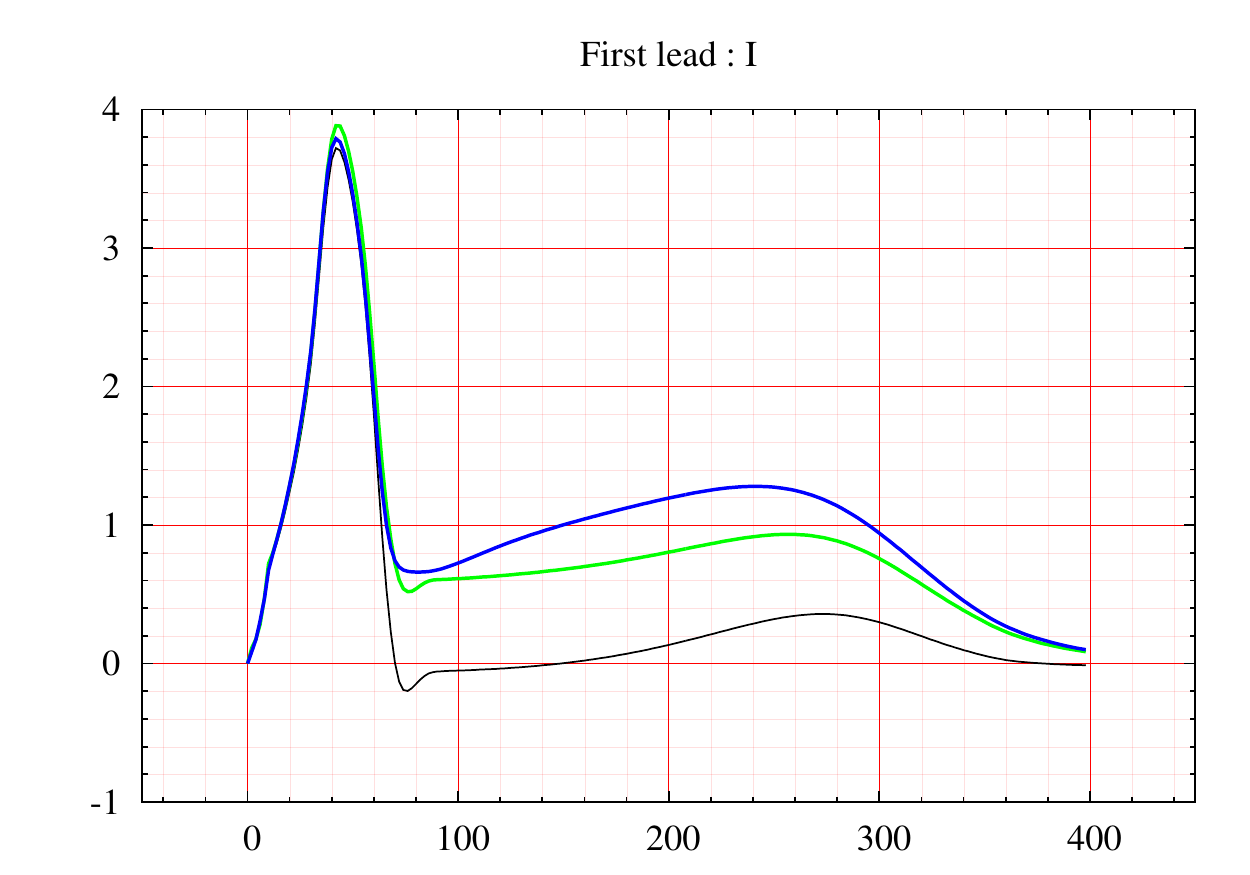}
                  \caption{First lead of the ECG for an infarction centered in point $P$: ECG with the complete model (green line) and ECG obtained using POD (blue line). The healthy reference case is given by the black line.}
                  \label{fig:ecgV1infarct}
                 \end{minipage}\\ \ \\ \ \\ 
\begin{minipage}{0.98\textwidth}
 \centering
 {
 \subfloat[t = 70ms]{ \label{fig:unknown60}
 \includegraphics[viewport = 0 0 1332 796 , width=0.42\textwidth]{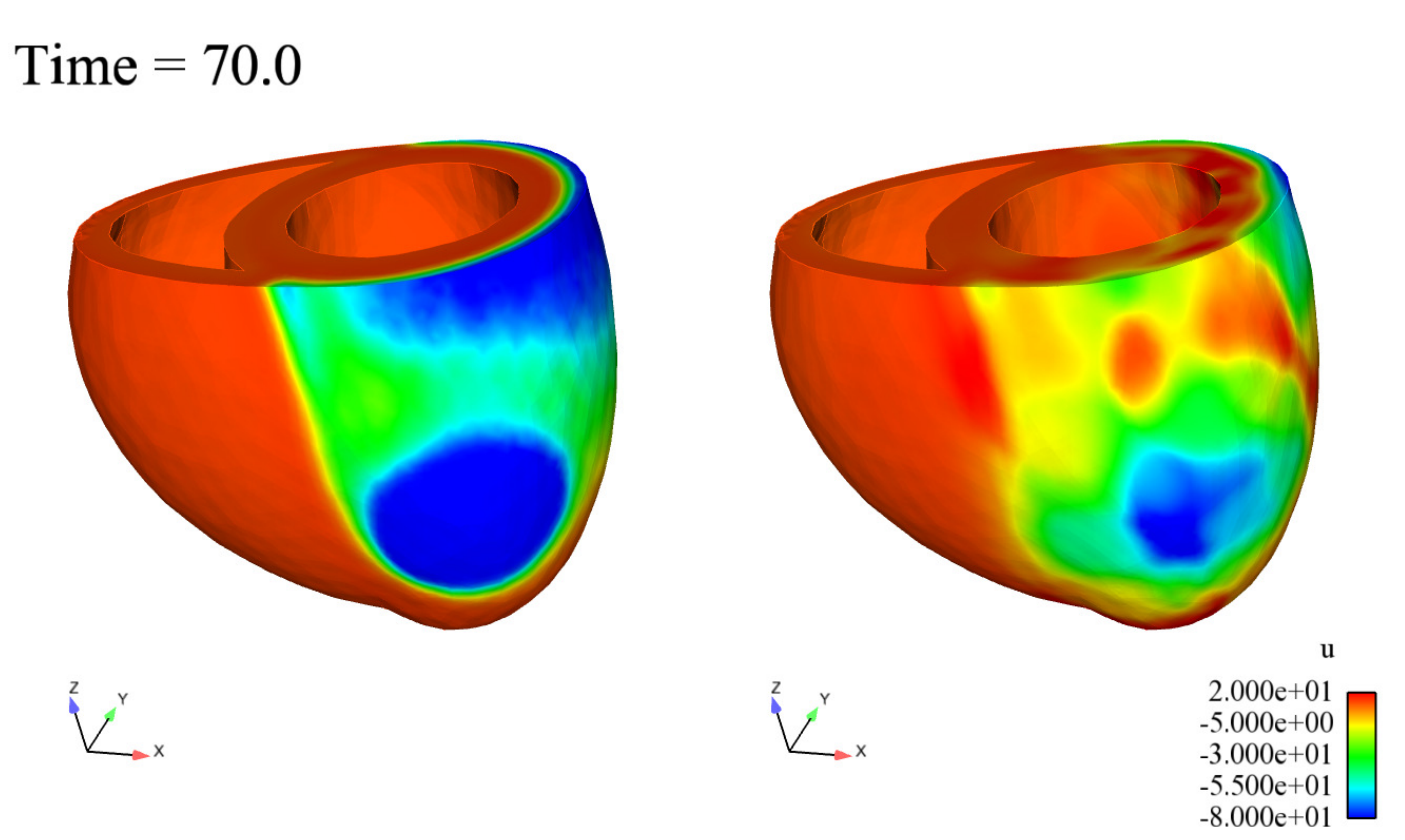}} \ \hspace{0.5 cm} \ 
 \subfloat[t = 200ms]{ \label{fig:unknown200}
 \includegraphics[viewport = 0 0 1332 796 , width=0.42\textwidth]{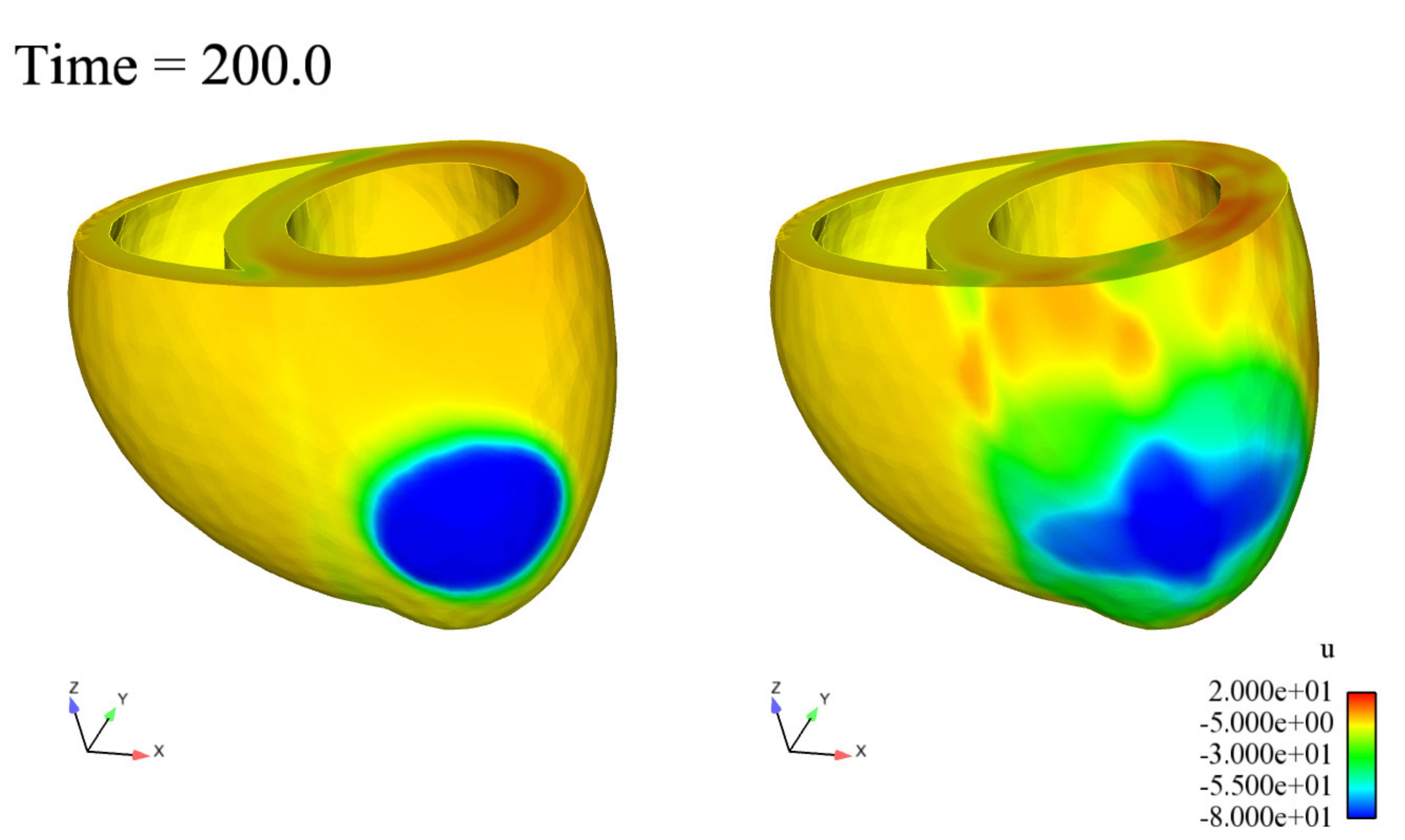}}
 }
 \caption{Simulation of an infarction in point $P$: using the complete model (left) and using the POD method (right).}
 \label{fig:infarct_unknown}
\end{minipage}
\end{figure}

We first generate the POD basis from a \emph{healthy} simulation : we run a 400 milliseconds simulation for an healthy test case, with a time-step of 0.5 milliseconds, keep snapshots every 2 milliseconds and obtain a basis of 100 vectors. Then, we use this basis to simulate an infarction centered in the arbitrary point $P$ indicated by the arrow in Figure \ref{fig:inf_points}.

\begin{figure}[htbp]
\begin{minipage}{1.\textwidth}
\centering
\includegraphics[width=0.9\textwidth]{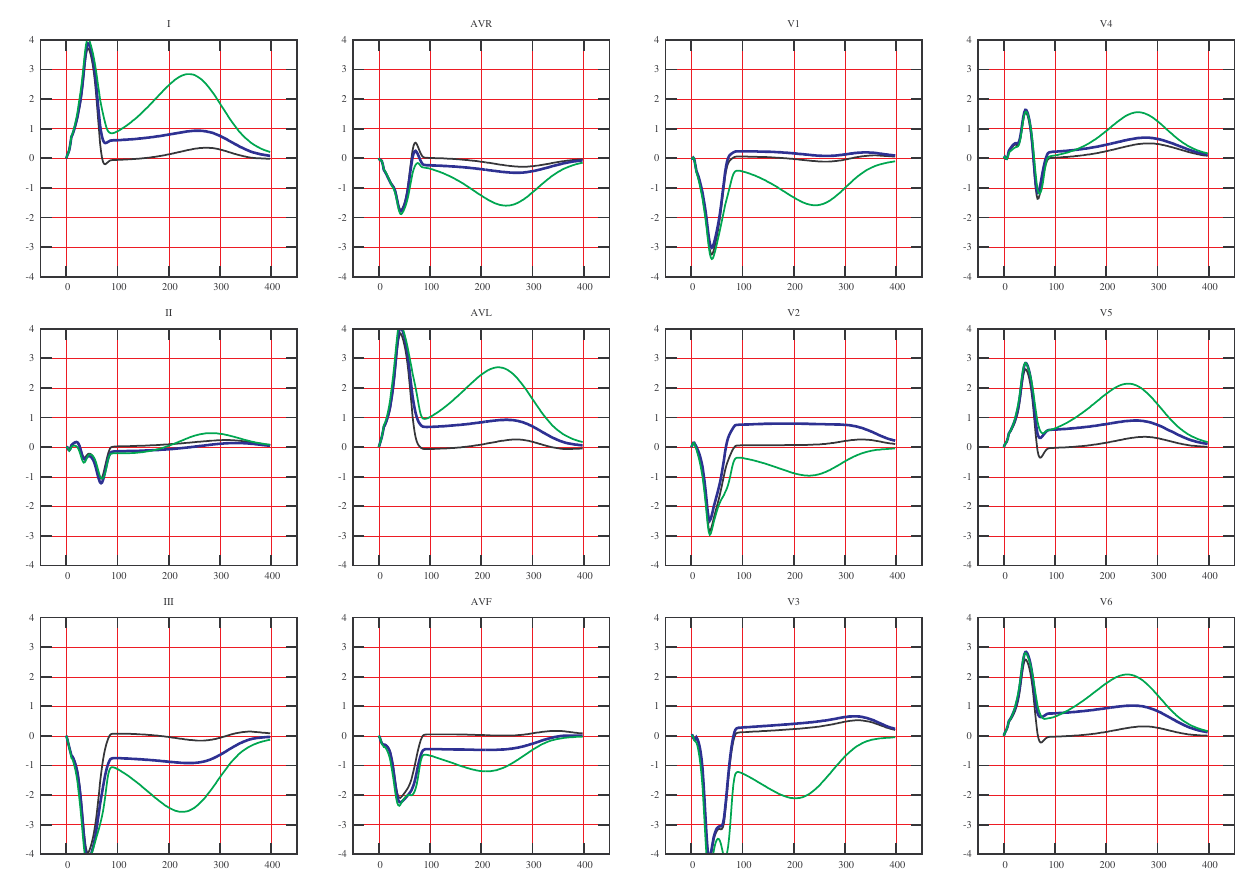}
  \caption{Simulated ECGs for an infarction centered in an arbitrary point $P$ with a POD basis generated from the healthy case: ECG with the full model (blue line), ECG with the POD (green line) and ``healthy'' full-order ECG  (black line).}
\label{fig:ecg_infarction_wrong} 
\end{minipage}\\
\centering
\begin{minipage}{0.8\textwidth}
 \centering
 {
 \subfloat[Complete model]{ \label{fig:ecg_infarction_wrong_complete}
 \includegraphics[viewport = 0 0 877 620 , width=0.42\textwidth]{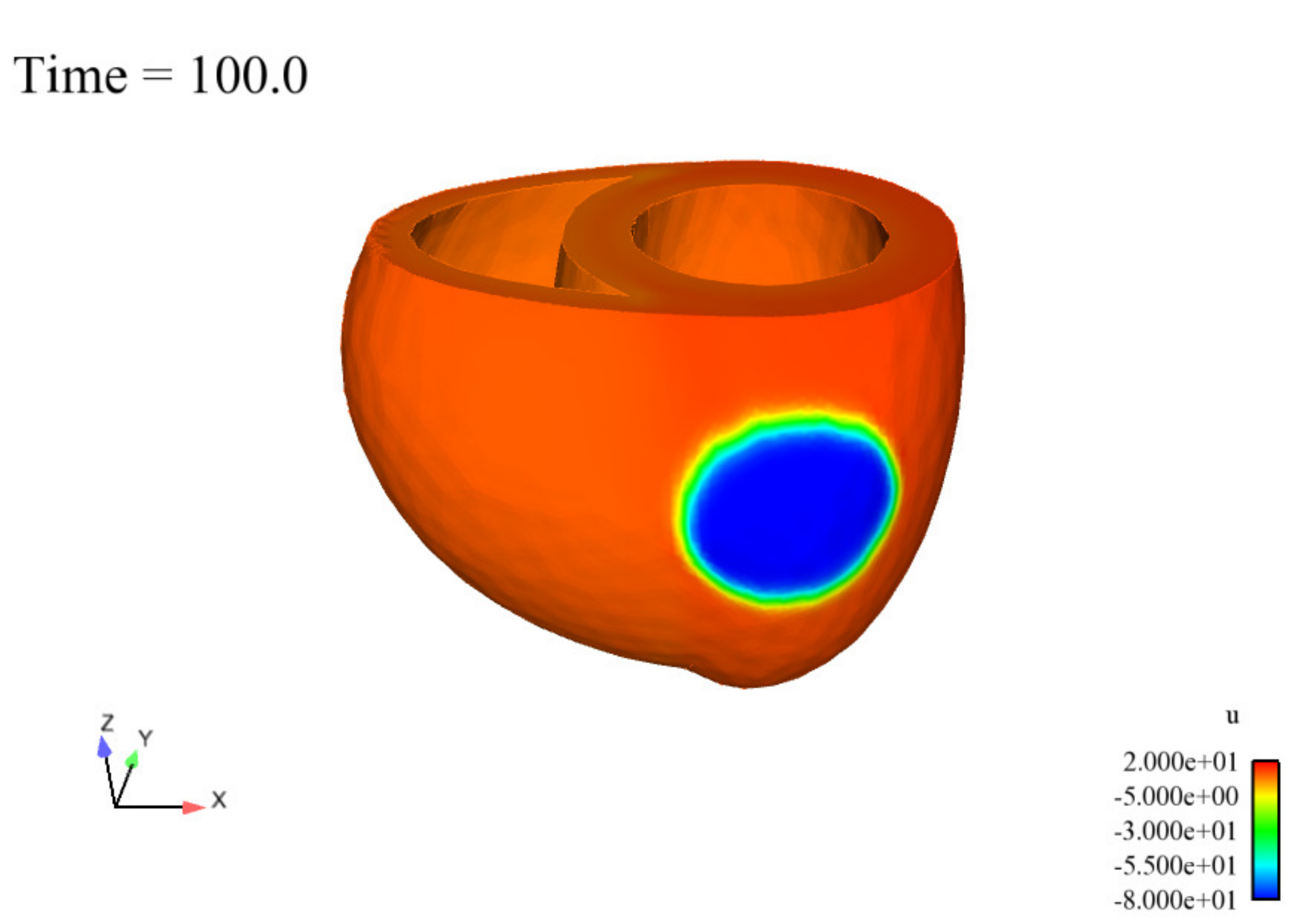}} \
 \subfloat[Reduced model]{ \label{fig:ecg_infarction_wrong_POD}
 \includegraphics[viewport = 0 0 877 620 , width=0.42\textwidth]{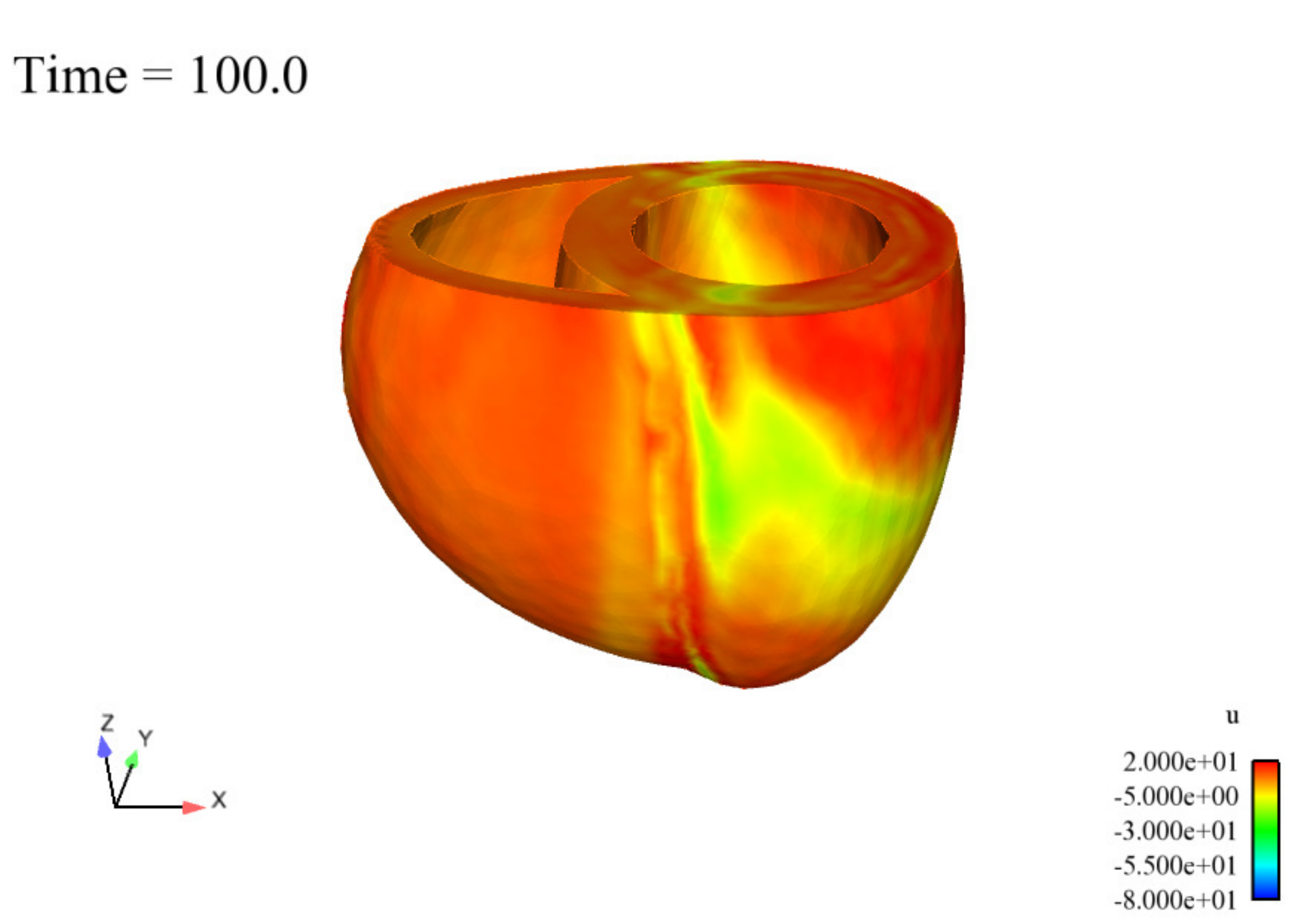}}
 }
 \caption{Simulation of anterior infarction using the complete model (left column) and the reduced model (right column) with a POD basis generated from the healthy case.}
 \label{fig:snapshots_infarction_wrong}
\end{minipage}
\end{figure}

A comparison of the green line (reduced-order model) and blue line (full-order model) in Figure~\ref{fig:ecg_infarction_wrong} shows that this basis does not allow to approximate accurately the ECG.

A look at the transmembrane potential (Figure \ref{fig:snapshots_infarction_wrong}) shows that the ``healthy'' POD basis is indeed unable to capture the sharp variation induced by the infarcted region. This explains the poor results observed on the ECG.

To improve the space approximation property of the POD basis, we propose to collect the snapshots for different infarction points. More precisely, 100 snapshots are taken from an healthy case and 50 snapshots for each infarction of the 18 points shown in Figure~\ref{fig:inf_points}. Moreover, since most of the solutions variation occurs during the first $100$ ms, half of the snapshots are taken in this period and the other half between $100$ ms and $400$ ms. Then a POD basis of 100 vectors is computed as usual. 

This new POD basis is used to simulate the infarction on the point $P$ indicated in Figure \ref{fig:inf_points}. A comparison of  the results obtained with the full and reduced-order models is given in Figure~\ref{fig:infarct_unknown}. Although the results look better than with a basis coming from a pure healthy case, we observe that the solution seems to superimpose the solutions coming from all the nearest infarction points. As a result, the infarcted area seems too large with the reduced-order model during the repolarization phase, which induces a difference of the ST elevations in the corresponding ECG. Indeed,  we see in Figure \ref{fig:ecgV1infarct} that the curves of the ECGs have a similar shape but different magnitude for the full and reduced-order models.  This discrepancy would probably be reduced by refining the grid of the precomputed infarcted regions.

\section{Application to inverse problems}\label{sec:inverse}

In this section, we propose a strategy to identify some parameters of the model from the ECG. The first application deals with the identification of the 4 parameters considered in Section~\ref{sec:pod-ex}. In the second application, we try to recover the location of an infarction modeled as in Section~\ref{sec:PODinfarct}. Note that in this preliminary study, the reference ECGs used for the identification are ``synthetic'' \emph{i.e.} previously computed by the model itself.

\subsection{Optimization method}

The parameter identification is done by minimizing the discrepancy between a reference ECG and the ECG provided by the model with a given set of parameters. 
More precisely, let us denote by $n\in \mathbb{N}^*$ the number of parameters and by $\theta \in  \mathbb{R}^n$ the vector of parameters we are looking for in a subset $\mathcal{I}$ of $\mathbb{R}^n$. The subset $\mathcal{I}$ is given by $\mathcal{I}_1\times\dots \times \mathcal{I}_n$ where $\mathcal{I}_j$ is an interval where the value $\theta_j$ is assumed to be. The following cost function is minimized 
\[
J(\theta)=\int_0^T|V_{\mathrm I}(t)-V_{\mathrm {I,ref}}(t)|^2+|V_{\mathrm {II}}(t)-V_{\mathrm {II,ref}}(t)|^2+|V_{\mathrm {III}}(t)-V_{\mathrm {III,ref}}(t)|^2\, dt
 \]
with respect to $\theta \in I$  where $V_{\mathrm I}$, $V_{\mathrm {II}}$ and $V_{\mathrm {III}}$ are the three Einthoven leads given by the simulation for the value $\theta$ of the parameters and $V_{\mathrm {I, ref}}$, $V_{\mathrm {II, ref}}$ and $V_{\mathrm {III,ref}}$ are the Einthoven leads of the reference ECG. 
This functional will be approximated by
\[
J(\theta)=\delta t \sum_{i=1}^N
\left[ 
|V_{\mathrm I,i}-V_{\mathrm {I,ref},i}|^2+|V_{\mathrm {II},i}-V_{\mathrm {II,ref},i}|^2+|V_{\mathrm {III},i}-V_{\mathrm {III,ref},i}|^2
\right]
 \]
where $W_i$ is the numerical approximation of $W(t_i)$ for $W=V_{\mathrm I}$, $V_{\mathrm {II}}$, $V_{\mathrm {III}}$, $V_{\mathrm {I,ref}}$, $V_{\mathrm{II,ref}}$ and $V_{\mathrm{III,ref}}$.
This cost function only uses information coming from the three Einthoven leads, numerical tests have been also run
with the twelve standard leads and give similar results (at least for the present test case).

This optimization problem is solved by an evolutionary algorithm (see \cite{dumas-elalaoui-07,dumas-2010}).
We briefly indicate its main steps, for the sake of completeness.
First, an initial population of elements, called individuals, is randomly created. Then, the algorithm modifies the population in order to promote the best individuals according to the cost function. Considering $N_p$ elements $(\theta_1,\dots,\theta_{N_p}) \in  I^{N_p}$ corresponding to different values of the parameters to identify, the  population is regenerated $N_g$ times, where $N_g$ corresponds
to  the number of generations. At each generation, $J$ is
evaluated  for each individual and the population evolves from
a generation to another following three stochastic principles inspired from the darwinian theory of evolution of species
 \begin{itemize}
 \item selection: promote the individuals whose value by the cost function is small, the members which cost function is smaller are preserved in the next generation, while those with a high cost function are killed;
 \item crossover: create from two individuals two new individuals by doing a barycentric combination of them with random and independent coefficients, a linear combination of the two selected vectors is done in order to create a new member $\theta \in  \mathbb{R}^n$;
 \item mutation: consist of replacing an individual by a new one randomly chosen in its neighborhood \emph{i.e.} a member whose values are closer to the referred one in the sense of an Euclidean norm. In our case, the amplitude of the mutation get smaller when the number of generations increases.
 \end{itemize}
At each generation, a one-elitism principle is added in order to make sure to keep the best individual of the previous population. 

Among its advantages, the genetic algorithm can easily be run in parallel. Moreover, unlike deterministic descent methods which require the computation of the gradient of the cost function, it can be very easily implemented. Its main flaw is to require a large number of evaluations of the direct problem, and since the initial population is chosen randomly, the minimization process  has to be run several times. To reduce the number of exact evaluations,  we have considered the approximate genetic algorithm based on a surrogate model strategy (see \cite{dumas-2010}). The idea is to approximate the cost function for a part of the population, taking advantage of the growing database of exact evaluations. The approximation is done by interpolation with Radial Basis Functions. The maximum number $N_{ex}$ of total exact evaluations is  fixed and the number of exact evaluations decreases at each generation. For instance, in the case of an initial population of $N_p=80$ members and a fixed number of $N_g=15$ generations, we impose the maximum total number of exact evaluations $N_{ex}=600$; typically we impose all evaluations to be exact during the first 4 or 5 generations, then the number of exact evaluations decreases constantly at each generation (see \cite{dumas-2010}).

Even with the approximate genetic algorithm, solving the optimization problem is still very time-consuming.  Our strategy is thus to speed up the evaluation of the cost function by using the reduced-order model. 

\subsection{Identification of four parameters}\label{sec:4param}

In this section, we test our identification strategy for the four parameters 
$\tau_{\mathrm {in}}$,   $C_{\mathrm {m}}$, $A_{\mathrm {m}}$ and $\tau_{\mathrm{close}}^{\mathrm{RV}}$.  The reference ECG used in the cost function is the one computed in  \cite{boulakia-cazeau-fernandez-gerbeau-zemzemi-10} obtained with the full-order model for the values  $(\tau_{\mathrm {in}},C_{\mathrm {m}},A_{\mathrm {m}},\tau_{\mathrm{close}}^{\mathrm{RV}})=(0.8,10^{-3},200,120)$. The values of $(\tau_{\mathrm {in}},C_{\mathrm {m}},A_{\mathrm {m}},\tau_{\mathrm{close}}^{\mathrm{RV}})$ are searched for in the set
$$
[0.5,1.5] \times [5\times 10^{-4},2\times 10^{-3}] \times [100,300] \times [50,150].
$$
The key point is to perform the ``exact'' evaluations required by the optimization algorithm with the reduced-order model. Results presented in Table \ref{tab:in-cm-am-close} are obtained with 80 POD modes and with the following parameters for the genetic algorithm:
$N_p=120$, $N_g=12$ and $N_{ex}=850$. Just to give an idea, the computational time was about 2 days on a PC with 16Go of RAM and using six processors Intel Xeon 3.2 GHz. Each one of the 800 time-steps of the ``exact'' evaluation requires about 3.5 seconds using the full-order model, while it is reduced to about 0.5 second using the reduced model. If the 850 ``exact'' evaluations were performed with the full-order model, the computational time will be unacceptable (more than two weeks). The time needed to compute the POD basis is negligible with respect to the time needed by the overall simulation.

As explained in section~\ref{sec:forward}, the accuracy of the reduced-order model is satisfactory when the four parameters of interest vary within a reasonable range. It is therefore possible to work with a POD basis generated for a single set of parameter $\theta_0$. This is the approach referred to as \textbf{M1} in Tables \ref{tab:in-cm-am-close} and \ref{tab:in-cm-am-close2}. 
Nevertheless, to improve the accuracy, it is possible to use many POD bases computed ``off-line'' for different values of $\theta$ taken in a finite subset $\mathcal{A}$ of $I$. Next, for a given value $\theta \in I$, the POD basis used for the resolution corresponds to the closest value of $\theta\in\mathcal{A}$. This is the approach \textbf{M2} in Tables \ref{tab:in-cm-am-close} and \ref{tab:in-cm-am-close2}. 

\begin{table}[htbp]
\centering
\[
\begin{array}{c|c|c|c}
& \hspace{0.2cm}(\tau_{\mathrm {in}},C_{\mathrm {m}},A_{\mathrm {m}},\tau_{\mathrm{close}}^{\mathrm{RV}})\hspace{0.2cm} & \textrm{Relative error
    (in \%)}& \textrm{Value of $J$}\\
\hline
\hspace{0.2cm}\textrm{{\bf M1} with } \theta_0=\theta_0^1 \hspace{0.2cm} & (0.95,9.3\times 10^{-4},185,126) &  9.6 & 3.05\\
\hline
\hspace{0.2cm}\textrm{{\bf M1} with } \theta_0=\theta_0^2 \hspace{0.2cm} & (0.93,1.05\times 10^{-3},162,128) &  11.7& 7.4\\
\hline
\hspace{0.2cm}\textrm{{\bf M2} with $\mathcal{A}$ given by (\ref{defA2})}  \hspace{0.2cm} & (0.86,10^{-3},179,123.5) & 5.2&2.15\\
\end{array}
\]
\caption{ Identification of $(\tau_{\mathrm {in}},C_{\mathrm {m}},A_{\mathrm {m}},\tau_{\mathrm{close}}^{\mathrm{RV}})$ (Reference value $(0.8,10^{-3},200,120)$).}
\label{tab:in-cm-am-close}
\end{table}
The relative error is the mean relative error given by the following formula:
\[
\frac{1}{4}\left(\frac{|\tau_{\mathrm {in}}-0.8|}{0.8}+\frac{|C_{\mathrm {m}}-10^{-3}|}{10^{-3}}+\frac{|A_{\mathrm {m}}-200|}{200}+\frac{|\tau_{\mathrm{close}}^{\mathrm{RV}}-120|}{120}\right).
\]
Here, $\theta_0^1= (1,10^{-3},200,100),\, \theta_0^2=(1.3,1.4\times 10^{-3},170,90)$ and in the last line of the table, we have considered the POD method {\bf M2} with $\mathcal{A}$ given by
\begin{equation}\label{defA2}	
\begin{array}{ll}
	
\mathcal{A}=&\big\{0.5;1;1.5 \} \times \{ 5\times 10^{-4};10^{-3};1.5\times 10^{-3};2\times 10^{-3}\big\} \\
& \times \big\{100;200;300 \big\}\times \big\{ 50;100;150\big\}.
\end{array}
\end{equation}

The genetic algorithm is run several times and the results shown corresponds to a mean value.
On this example, we see that method {\bf M2} allows approximately to divide the error  by $2$. The price of this improvement was to compute $108$ POD bases. Although this computation was done off-line, this approach therefore requires a significant computational effort. A larger population has also been considered in Table \ref{tab:in-cm-am-close2} with $N_p=300$, $N_g=20$ and $N_{ex}=1700$. Here again, we notice that method {\bf M2} allows to improve the identification results.

\begin{table}[htbp]
\centering
\[
\begin{array}{c|c|c|c}
& \hspace{0.2cm}(\tau_{\mathrm {in}},C_{\mathrm {m}},A_{\mathrm {m}},\tau_{\mathrm{close}}^{\mathrm{RV}})\hspace{0.2cm} & \textrm{Relative error
    (in \%)}& \textrm{Value of $J$}\\
\hline
\hspace{0.2cm}\textrm{{\bf M1} with } \theta_0=\theta_0^1 \hspace{0.2cm} & (0.83,1.02\times 10^{-3},184.2,123.1) & 4.1  & 2.2\\
\hline
\hspace{0.2cm}\textrm{{\bf M1} with } \theta_0=\theta_0^2 \hspace{0.2cm} & (0.91,1.09\times 10^{-3}, 153,126.4) &  12.3& 7.55\\
\hline
\hspace{0.2cm}\textrm{{\bf M2} with $\mathcal{A}$ given by (\ref{defA2})}  \hspace{0.2cm} & (0.83,1.01\times 10^{-3},189,123.2) & 3&2.05\\
\end{array}
\]
\caption{ Identification of $(\tau_{\mathrm {in}},C_{\mathrm {m}},A_{\mathrm {m}},\tau_{\mathrm{close}}^{\mathrm{RV}})$ (Reference value $(0.8,10^{-3},200,120)$) with a larger population.}
\label{tab:in-cm-am-close2}
\end{table}

\subsection{Identification of an infarcted zone}\label{sec:identinf}

We are interested in estimating the location of an infarcted area, modeled as explained in section~\ref{sec:PODinfarct}. We proceed as before by minimizing the discrepancy between a reference and a simulated ECG, using the genetic algorithm. 

\begin{figure}[!htp]
\begin{minipage}{0.9\textwidth}
\centering
\includegraphics[width=1.\textwidth]{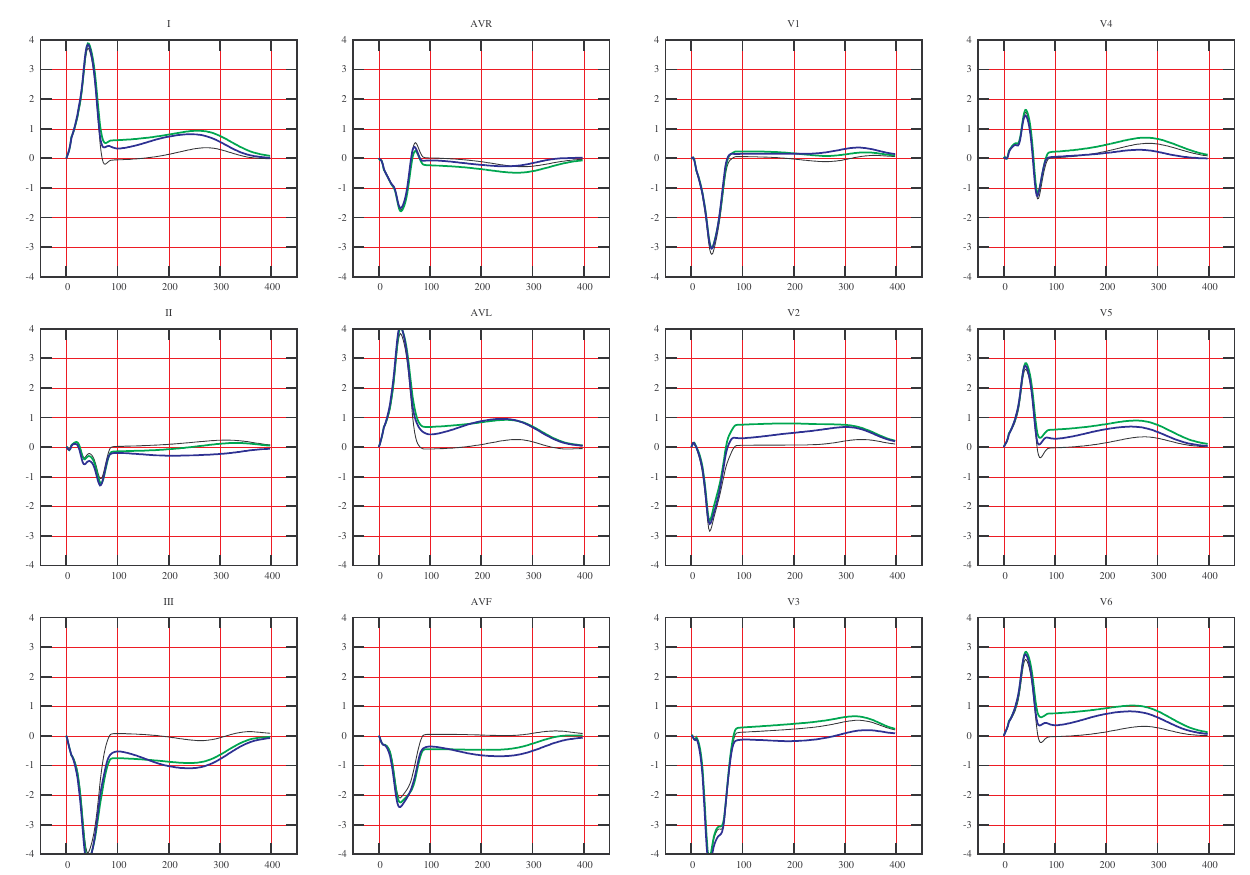}
\caption{Simulated ECG for an infarction located in point $P$: green line represents the simulated reference ECG and blue line the ECG corresponding to the infarcted center found with the resolution of a genetic algorithm. Black line gives the healthy reference case.}
\label{fig:ecg_inf_identification}
\end{minipage}\\ \ \\
\begin{minipage}{0.98\textwidth}
 \centering
 {
 \subfloat[t = 80ms]{ \label{fig:inf_identification60}
 \includegraphics[viewport = 0 0 877 620 , width=0.4\textwidth]{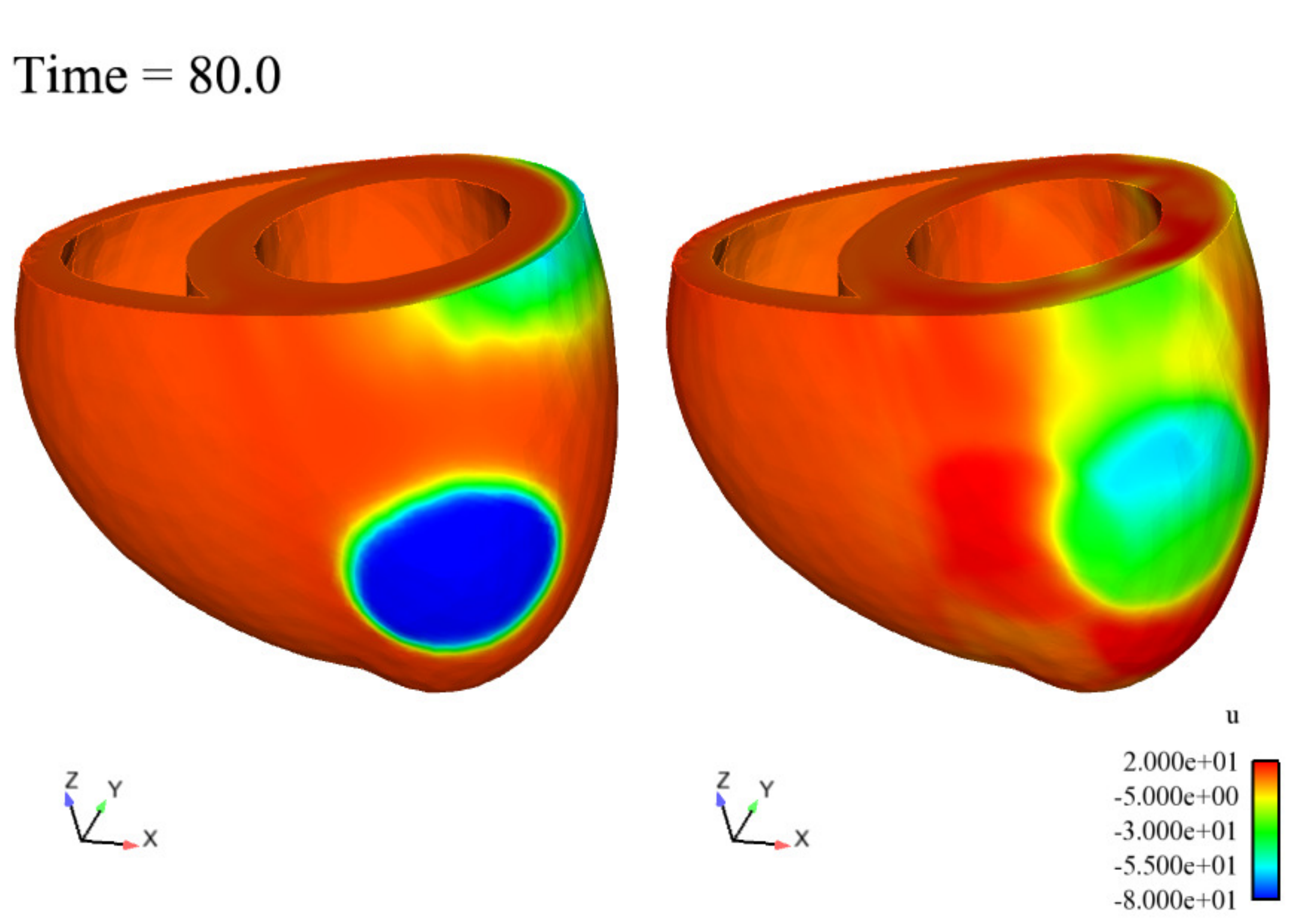}} \hspace{0.3 cm}\\
 \subfloat[t = 300ms]{ \label{fig:inf_identification200}
 \includegraphics[viewport = 0 0 877 620 , width=0.4\textwidth]{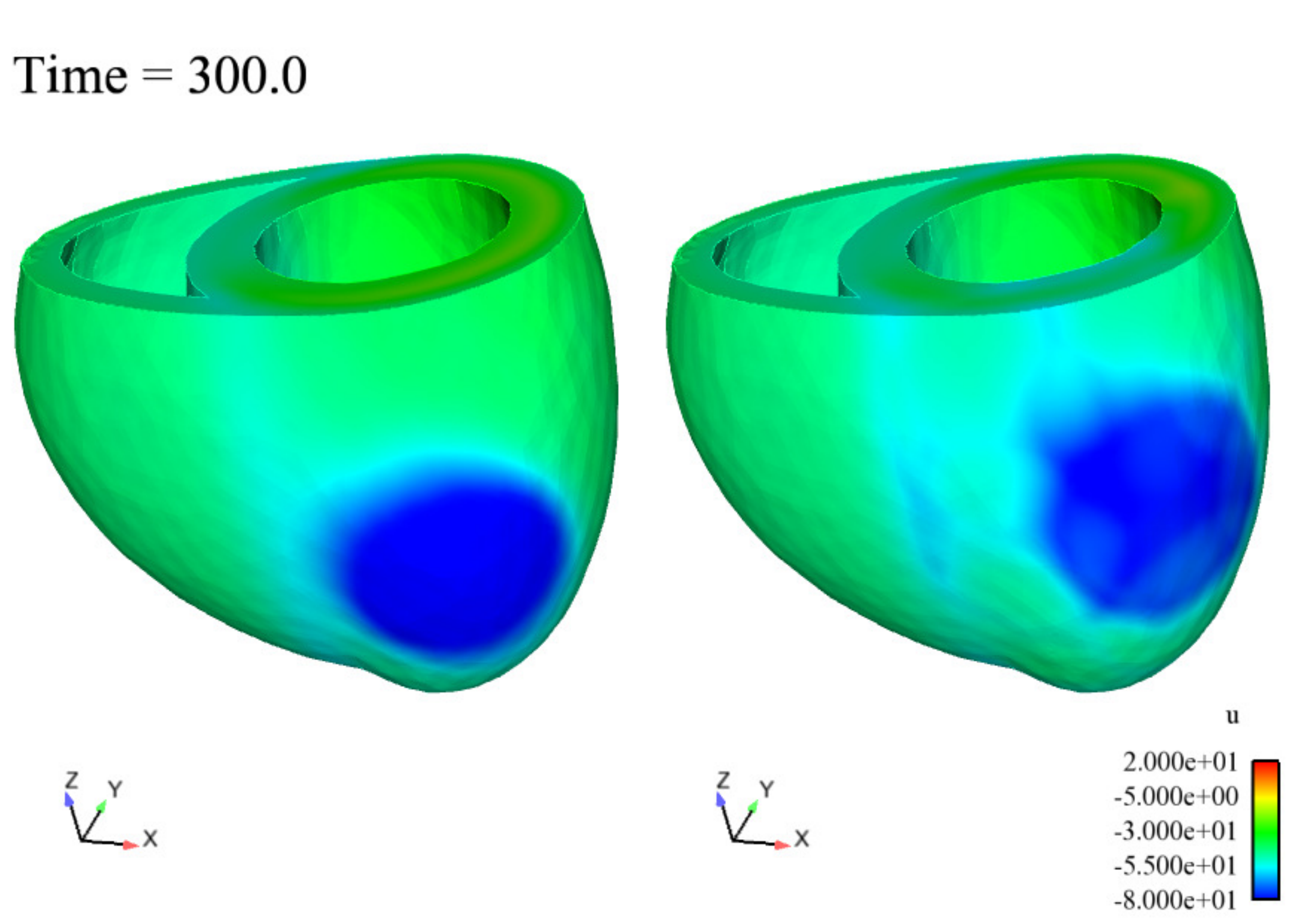}}
 }
 \caption{Left: transmembrane potential calculated in the reference case solved with the full model. Right: transmembrane potential of the solution found with the genetic algorithm obtained with the POD.}
 \label{fig:inf_identification}
\end{minipage}
\end{figure}

The reference ECG is obtained solving the complete model for an infarction area centered at point $P$ as described in section \ref{sec:PODinfarct}. The genetic algorithm is run with $N_p=80$, $N_g=15$ and $N_{ex}=600$, and the parameters we try to identified are $\theta=(x_P,y_P,z_P)$. The point $P$ is only searched for in the left ventricle domain.

The results are reported in Figures \ref{fig:ecg_inf_identification} and \ref{fig:inf_identification}. The reference ECG (blue line in Figure \ref{fig:ecg_inf_identification}) is well approximated by the one obtained from the resolution of the genetic algorithm (blue line). The identified infarcted region is actually very close to the reference one (Figure \ref{fig:inf_identification}). The solution of the genetic algorithm can be improved by including more off-line experiments, as indicated at the end of section \ref{sec:PODinfarct}.

For a different approach to tackle this problem, we refer to~\cite{nielsen-lysaker-tveito-07}. An interesting possibility to investigate could be to enrich the cost function giving more weight to the ST deviation, as suggested in \cite{goletsis-papaloukas-fotiadis-likas-michalis-04}.


\section{Conclusion}\label{sec:conclusion}

We have presented some preliminary results obtained with a reduced-order model of electrophysiology based on the Proper Orthogonal Decomposition (POD) method. A well-known difficulty of reduced-order modeling is to identify those situations where the reduced-model can be considered as reliable. We do not claim that we have answered this difficult question in this paper. Nevertheless, our experiments may suggest some general trends. First, the reduced-order model seems quite robust when the ionic current parameters are perturbed and when the heart beat increases. This may have interesting applications for long time simulations, for example to compute a restitution curve, or to estimate the ionic current parameters in an optimization loop. Second, a reduced-order model based on a single simulation is totally inadequate to approximate a situation with a \emph{spatial} change in the parameter: this has been shown in the present study for an infarcted region; in~\cite{boulakia-gerbeau-11} the same observation was made for the initial activation region. A natural strategy consisting of precomputing several POD bases with different sets of parameters, or a single POD basis from different experiments, proved to be satisfactory in the cases we have considered. Nevertheless, this solution requires an important off-line effort which makes it difficult to apply with more than a few parameters. Alternative strategies have therefore to be investigated in future works.

\bibliographystyle{plain} 
\bibliography{biblio}

\tableofcontents

\end{document}